\documentclass[11pt]{article}
\usepackage{amsfonts}
\usepackage{epsfig,graphics,subfigure,psfrag,amsmath,amssymb}
\usepackage{latexsym,wrapfig}
\usepackage{bm}
\usepackage{amsthm}
\usepackage{graphicx}
\usepackage{mathrsfs}
\usepackage{dcolumn}
\usepackage{color}
\usepackage{pstricks}
\usepackage{epstopdf}
\usepackage{subfigure}
\usepackage{enumerate}
\usepackage{caption}
\usepackage{geometry}
\usepackage{threeparttable}

\usepackage{hyperref}
\makeatletter
\@addtoreset{equation}{section}
\makeatother

\usepackage{latexsym,amsmath,amssymb,paralist}

\usepackage{comment}
\usepackage{algorithm}
\usepackage{algpseudocode}
\usepackage[title]{appendix}
\usepackage{graphicx}
\usepackage{epstopdf}
\usepackage{bm}% bold math
\newcommand{\reff}[1]{{\rm (\ref{#1})}}

\newcommand{\R}{\mathbb{R}}            % real numbers
\newtheorem{thm}{Theorem}[section]

\newtheorem{lem}[thm]{Lemma}

\newtheorem{defn}[thm]{Definition}

\newtheorem{rem}[thm]{Remark}
\def\E\mathbb{ E}

\newcommand{\be}{\begin{eqnarray}}
\newcommand{\ee}{\end{eqnarray}}
\newcommand{\ben}{\begin{eqnarray*}}
\newcommand{\een}{\end{eqnarray*}}

\def\XXint#1#2#3{{\setbox0=\hbox{$#1{#2#3}{\int}$}
\vcenter{\hbox{$#2#3$}}\kern-.51\wd0}}

{\allowdisplaybreaks
\begin{document}

\title{Structure-Preserving Numerical Methods for Nonlinear Fokker--Planck Equations with Nonlocal Interactions by an Energetic Variational Approach}

\author{Chenghua Duan\footnotemark[1]\and Wenbin Chen\footnotemark[2]\and Chun Liu\footnotemark[3]\and Xingye Yue\footnotemark[4] \and Shenggao Zhou\footnotemark[5]}
% .

 \renewcommand{\thefootnote}{\fnsymbol{footnote}}

 \footnotetext[1]{Shanghai Center for Mathematical Sciences,  Fudan University, Shanghai 200438, China. Email:  chduan@fudan.edu.cn.}
  \footnotetext[2]{School of Mathematical Sciences and Shanghai Key Laboratory for Contemporary Applied Mathematics, Fudan University, Shanghai 200433, China. Email: wbchen@fudan.edu.cn.}
 \footnotetext[3]{Department of Applied Mathematics, Illinois Institute of Technology, Chicago, IL 60616, USA. Email: cliu124@iit.edu.}
 \footnotetext[4]{Department of Mathematics and Mathematical Center for Interdiscipline Research, Soochow University, Suzhou 215006, Jiangsu, China. Email: xyyue@suda.edu.cn.}
 \footnotetext[5]{Corresponding author. Department of Mathematics and Mathematical Center for Interdiscipline Research, Soochow University, Suzhou 215006, Jiangsu, China. Email: sgzhou@suda.edu.cn.}
\date{\today}
\maketitle

\begin{abstract}
In this work, we develop novel structure-preserving numerical schemes for a class of nonlinear Fokker--Planck equations with nonlocal interactions. Such equations can cover many cases of importance, such as porous medium equations with external potentials, optimal transport problems, and aggregation-diffusion models.  Based on the Energetic Variational Approach,  a trajectory equation is first derived by using the balance between the maximal dissipation principle and least action principle.   By a convex-splitting technique, we propose energy dissipating numerical schemes for the trajectory equation.   Rigorous numerical analysis reveals that the nonlinear numerical schemes are uniquely solvable,  naturally respect mass conservation and positivity at fully discrete level, and preserve steady states.  Under certain smoothness assumptions, the numerical schemes are shown to be second order accurate in space and first order accurate in time. Extensive numerical simulations are performed to demonstrate several valuable features of the proposed schemes.  In addition to the preservation of physical structures, such as positivity, mass conservation, discrete energy dissipation, blue and steady states, numerical simulations further reveal that our numerical schemes are capable of solving \emph{degenerate} cases of the Fokker--Planck equations effectively and robustly. It is shown that the developed numerical schemes have convergence order even in degenerate cases with the presence of solutions having compact support, can accurately and robustly compute the waiting time of free boundaries without any oscillation, and can approximate blow-up singularity up to machine precision.

%The solution of the system is expected to satisfy three features: positivity, the conservation of mass and the energy dissipation law. In this paper, based on an energetic variational approach, a force balance between the least action principle and the maximum dissipation principle,  the trajectory equation of the system can be obtained. Then by the convex splitting scheme, we can establish a numerical scheme that can preserve the positivity of the solution and the conservation of mass naturally. Moreover, we can prove that the scheme is unique solvable on a convex set and preserves the corresponding discrete energy dissipation law. Meanwhile, under some smooth assumption, it is proved that the scheme is second-order convergent  in space and  first-order convergent  in time by a higher order expansion technique. The numerical results not only confirm the order for smooth test, but also show the order for compact test.  The numerical solutions preserve long time asymptotics effectively. Furthermore, there are other advantages: dealing  with the nonlinear mobility naturally, catching the blow-up solution at a machine precise over any equidistant grid and predicting the possible waiting time phenomenon.

\bigskip


\end{abstract}

\section{Introduction}
We focus on the following initial-boundary value problem
\begin{equation}\label{equ:GFP}
\left\{
\begin{aligned}
   & \partial_t u=\partial_x\{f(u)\partial_x[H'(u)+V(x)+W*u]\}, \ x\in\Omega,\ t>0,\\
   & u(x,0)=u_0(x),\ x\in\Omega, \\
   & f(u)\partial_x[H'(u)+V(x)+W*u]=0, \ x\in\partial\Omega, \ t>0,
\end{aligned}
\right.
\end{equation}
where $u(x,t)\geq 0$ represents the time-dependent probability density,
$\Omega \subset \mathbb{R} $ is a bounded domain,
$H(\cdot): \R^{+} \cup \{0\} \rightarrow \R$ is the density of internal energy with $H''(\cdot)>0$, $V(\cdot)$ is an external potential,  $W(\cdot)$ is an even Lipschitz continuous function describing particle interactions,
and $f: \R^{+} \cup \{0\} \rightarrow \R^{+} \cup \{0\}$ is a given increasing differentiable function with $f(0)=0$ and $f'(0)\neq 0$.
% In addition, we assume that $f(u)$ increases with respect to $u$ and satisfies $\frac{f(u)}{u}\leq C$ for some fixed positive constant $C$ \cite{Z.Sun(2018)}.

The Fokker-Planck (FP) equation in the problem \eqref{equ:GFP} arises from various applications. Such an equation can be derived as mean-field limits of particle systems and has been used in various models to describe interacting gases \cite{J.Carrillo(2003),C.Villani(2003)}, granular materials \cite{D.Benedetto(1998)}, collective motion of animals \cite{C.M.Topaz(2006),T.Kolokolnikov(2013),J.Carrillo(2010)}, and cell migration and chemotaxis phenomena in biology \cite{J.Carrillo(2018),P.M.Lushnikov(2008),M.Burger(2007)}.  The FP equation covers many cases of importance. For instance, when $f(u)=u$, $V=W=0$, and $H'(u)=u^m (m>1)$, it becomes the porous medium equation \cite{J.Carrillo(2000)}.
%%The equation  \eqref{equ:GFP} is a generalized formation arising from two typical scenarios.
%For instance, when $W=0$ and $f(u)$ is nonlinear, it corresponds to the nonlinear (possibly degenerate) Fokker--Planck (FP) equation,
%%\begin{equation}\label{equ:NDP}
%%\partial_t u=\partial_x\{f(u)\partial_x[H'(u)+V(x)]\},
%%\end{equation}
%which is related to many classical problems, such as the heat equation and porous medium equation.
When considering a nonzero nonlocal interaction term, i.e., $W\neq 0$,  it is referred as the nonlocal FP equation \cite{J.Carrillo(2015), H.Liu(2016), Y.Qian(2019)}.
%\begin{equation}\label{equ:NN}
%\partial_t u=\partial_x\{u\partial_x[H'(u)+V(x)+W*u]\}.
%\end{equation}
Typical interaction potentials, $W(\cdot)$, appearing in above applications %The nonlinear diffusion type models become  interesting problem recently\cite{M. Burger(2014),J. Carrillo(2015),J. Carrillo(2018)}.
include fully attractive cases,  such as the Newtonian or Bessel potentials in chemotaxis
\cite{J.Carrillo(2018)} and power-laws in granular materials \cite{J.Carrillo(2003)}; cases that are repulsive in the short range and attractive in the long range,
such as combinations of power-law potentials and Morse-type potentials in swarming \cite{J.Carrillo(2014),C.M.Topaz(2006)}; and cases with compactly supported potentials in many biological applications, such as networks and cell sorting \cite{J.Barre(2017),M.Burger(2014),J.Carrillo(2018)}.
%\cite{J.A.Carrillo(2018),J.Barre(2017),M.Burger(2014),J.Carrillo(2018)}.

Mathematically, any solution to the  problem~\eqref{equ:GFP} has three main properties
\begin{itemize}
\item Non-negativity: if $u_{0}(x)\geq0$, then $u(x,t)\geq0$, $\forall x \in \Omega$, $t>0$;
\item Mass conservation: $\int_{\Omega}u_{0}(x)dx= \int_{\Omega}u(x,t)dx$;
\item Energy dissipation:  \begin{equation}\label{eqEne}\frac{d}{dt} E^{total}=-\Delta\leq 0,\end{equation}
where
\[
E^{total}:=\int_{\Omega}H(u(x))dx+\int_{\Omega}u(x)V(x)dx
+\frac{1}{2}\int_{\Omega}\int_{\Omega}W(x-y)u(x)u(y)dydx,
\]
and $$\Delta=\int_{\Omega}f(u)\big|\partial_x [H'(u)+V(x)+W*u]\big|^2dx.$$
\end{itemize}

The property \eqref{eqEne} has played a critical role in analyzing the dynamics of the problem~\eqref{equ:GFP}  in the works \cite{J.Carrillo(2001),J.Carrillo(2003),C.Villani(2003),Z.Sun(2018)}. Therefore, it is crucial and highly desirable to develop numerical methods that are able to maintain an analogous energy dissipation in the discrete sense.  Another challenge to obtain physically faithful numerical solutions lies in the development of numerical schemes that can guarantee the non-negativity of the numerical density while retaining the mass conservation, especially in the degenerate case.  Recently, various numerical schemes addressing above concerns, ranging from finite volume methods to discontinuous Galerkin (DG) methods,  have been developed to numerically solve the FP type of equations in the literature. Finite volume schemes with second order accuracy have been proposed for the problem~\eqref{equ:GFP} in \cite{M.Bessemoulin-Chatard(2012), J.Carrillo(2015)}. The schemes have semi-discrete (in space) entropy dissipation and positivity preserving properties for explicit-in-time discretization under a restriction on time step size due to the Courant--Friedrichs--Lewy (CFL) condition.  Entropic schemes have been developed in \cite{C.Buet(2004)} to solve the FP equations for a simplified model of granular media. It has been proved that the entropic schemes have many attractive properties, such as mass conservation, entropy decay, and positivity and equilibrium preserving. Based on entropic average fluxes,  another type of entropic schemes have been constructed to solve the nonlocal, nonlinear FP equations \cite{L.Pareschi(2007)}. It also has been shown that the constructed entropic schemes are able to preserve positivity, semi-discrete entropy dissipation, and asymptotic steady states with arbitrary accuracy. The work \cite{H.Liu(2016)} has proposed high order direct DG schemes, in which a discrete version of entropy dissipation law is respected by numerical solutions and positivity is enforced by a delicate reconstruction algorithm that is able to maintain accuracy.   To achieve high order accuracy, high order DG schemes for~\eqref{equ:GFP} have been established in the work \cite{Z.Sun(2018)}. For an interaction potential with a smooth kernel, the proposed semi-discrete DG scheme admits an entropy inequality at discrete level. The fully discretized DG scheme is able to produce non-negative solutions under a time step size constraint, with the help of a positivity-preserving limiter.  Based on harmonic-mean approximations,  finite difference schemes that are proved to respect mass conservation and unconditional positivity preservation have been proposed in \cite{Y.Qian(2019)}. Estimates on the condition number of the coefficient matrix has been established as well.   More recently, a fully discrete, implicit-in-time finite volume scheme that ensures the positivity and energy-decaying properties has been established in the work \cite{R.Bailo(2018)}.

Another closely related model, the Poisson--Nernst--Planck (PNP) equations, can be regarded as the problem \reff{equ:GFP} with $f(u) = u$, $H(u) = u\log u$, and a nonlocal Coulombic interaction kernel that is coupled through a Poisson equation. Related numerical methods \cite{LiuWangJCP(2014), MettiXuLiuJCP(2016), LiuWangJCP(2017), DingWangZhou_JCP19} with structure-preserving properties for the PNP equations can be extended to numerically solve the the problem \reff{equ:GFP} as well. There are other types of numerical methods for the nonlocal case, e.g., particle methods \cite{J.Carrillo(2017)} and evolving diffeomorphisms methods \cite{J.Carrillo(2016)}.
% a mixed finite element method \cite{M.Burger(2010)}, and a blob method \cite{K.Craig(2016)}.
However, energy dissipation law in fully discrete level and convergence order of numerical schemes have not been well studied.  In addition to above structure-preserving features at fully discrete level, it is rather challenging to develop numerical schemes that can capture finite-speed propagation and possible waiting time in degenerate cases. For solutions with compact support, it is non-trivial to show the convergence order, even numerically. When a solution blows up in a finite time,  the standard finite difference methods,  finite volume schemes, or DG methods, only present the order of $\mathcal{O}(1/h)$ blow-up on an equidistant mesh with grid spacing $h$. Improvement addressing these issues is still in lack.  %with small positive $\varepsilon$ close to the machine precision.
% If the interaction potential is not involved, or the interaction is defined by a smooth kernel,

In this paper, we propose a novel numerical scheme based on an Energetic Variational Approach (EnVarA), which is a balance between the maximal dissipation principle (MDP) and least action principle (LAP). The approach was originated from a pioneering work due to Onsager \cite{L.Onsager(1931), L.Onsager1(1931)} and further improved by Strutt \cite{J.W.Strutt(1873)}.  In recent years, it has been applied to develop mathematical models for complex physical systems \cite{C.Liu(2003), Y.Hyon(2010), Q.Du(2009), B.Eisenberg(2010)}, as well as numerical schemes for porous medium equations \cite{C.H.Duan(2018)} and the Wright-Fisher model that describes genetic drift \cite{C.H.Duan(2017)}. We first derive a trajectory equation and then establish its numerical scheme by a convex splitting technique. The positivity and mass conservation of the numerical solution can be preserved naturally. Numerical analysis proves that the numerical scheme is uniquely solvable, satisfies a discretized energy dissipation law, and preserves steady states. The proposed numerical scheme for the trajectory equation can also be justified at theoretical level that the convergence rate is first order in time and second in space.  We conduct extensive numerical tests to demonstrate several valuable advantages of the proposed schemes in overcoming the difficulties in the development of numerical methods for the FP equations. In addition to the success in  preservation of physical structures, including positivity, mass conservation, discrete energy dissipation, and steady states, our numerical simulations further demonstrate that the proposed schemes are able to solve the \emph{degenerate} FP equations effectively and robustly. Numerical results reveal that the developed numerical schemes have convergence order even in degenerate cases with the presence of solutions having compact support, can accurately and robustly calculate the waiting time of free boundaries without any oscillation, and can simulate blow-up singularity up to machine precision.
%To further demonstrate the robustness, our numerical methods are applied to study various cases with compactly supported nonlocal kernels.
%Numerical approximation of the density becomes first order in space due to the boundary treatment that is designed to preserve mass conservation.
%Moreover, for the compacted supported nonlocal kernel, it has been shown that  the number of bumps of the steady state  and its relationship with the initial data and the diffusion coefficient. Furthermore, the convergence order for solutions with a compact support is also given numerically,
%For the smooth kernel, the numerical ysolutions show the  process and verify the analysis results.

This paper is organized as follows. The EnVarA and trajectory equation of the nonlinear Fokker--Planck equations are outlined in Section \ref{sec:2}. The numerical scheme is described in Section \ref{sec:3}. Subsequently, the proof of unique solvability, energy stability,  optimal rate convergence analysis and steady-state preserving is provided in Section \ref{sec:4}. Section \ref{sec:5} presents various numerical results. Finally, in Section \ref{sec:6}, we draw conclusions.

\section{Energetic Variational Approach}
\label{sec:2}
%In this section, we first derive the trajectory equation of \eqref{equ:GFP} by an energetic variational approach  (EnVarA).
%In the Wright-Fisher model, $x\in[0,1]$ and $f(x,t)\geq0$ can be viewed as the position of particles and the density of $x$ at time $t$, respectively. We first introduce  the different coordinate systems.

We first introduce the Lagrangian and Eulerian coordinate systems.
 \begin{defn}
 Suppose that $\Omega_{0}^{X}$ and $\Omega_{t}^{x}$ $\subset\mathbb{R}^{m}$, $m\in\mathbb{N}^{+}$, are domains with smooth boundaries, time $t>0$, and ${\bf v}$ is a smooth vector field in $\mathbb{R}^{m}$. The flow map $x(X,t):\Omega_{0}^{X}\rightarrow\Omega_{t}^{x}$ is defined as a solution of
 \begin{equation}\label{equ:v}
 \left\{
 \begin{aligned}
   & \frac{d}{dt}x(X,t)={\bf v}(x(X,t),t),\ \ t>0 , \\
   & x(X,0)=X ,
 \end{aligned}
 \right.\
 \end{equation}
 where $X=(X_{1},...,X_{m})\in\Omega_{0}^{X}$ and $x=(x_{1},...,x_{m})\in\Omega_{t}^{x}$. The coordinate system $X$ is called the Lagrangian coordinate and $\Omega_0^X$ is called the reference configuration;  the coordinate system $x$ is called the Eulerian coordinate and $\Omega_{t}^{x}$ is called the deformed configuration.
 % and $\textbf{u}$ is the Eulerian velocity.
 \end{defn}
Since $\Omega_{0}^{X}$ and $\Omega_{t}^{x}$ are the same domain described by different coordinate systems, we denote the domain under consideration uniformly by $\Omega$ in the rest of this paper. Also, we assume that the Jacobian of the flow map, $\det\frac{\partial x(X,t)}{\partial X}$, remains positive in time evolution.

Now we derive a trajectory equation for the nonlinear nonlocal Fokker--Planck equations. The initial-boundary value problem~\eqref{equ:GFP} is equivalent to
%\begin{equation}\label{equ:GFP-1}
%\left\{
%\begin{aligned}
%   & \partial_t u=\partial_x\{f(u)\partial_x[H'(u)+V(x)+W*u]\}, \ x\in\Omega,\ t>0,\\
%   & u(x,0)=u_0(x),\ x\in\Omega, \\
%   & f(u)\partial_x[H'(u)+V(x)+W*u]=0, \ x\in\partial\Omega, \ t>0,
%\end{aligned}
%\right.
%\end{equation}
\begin{equation}\label{equ:GFP-1}
\left\{
\begin{aligned}
 & \partial_t u+ \partial_x(u\textbf{v})=0 , \ x\in\Omega, \ t>0,\\
 &  \textbf{v}=-\frac{f(u)}{u}\partial_x[H'(u)+V(x)+W*u], \ x\in\Omega,   \\
 & u(x,0)=u_{0}(x)>0,\ x\in \Omega ,  \\
 & \left. f(u)\partial_{x}\{[H'(u)+V(x)+W*u]\}\right |_{\partial\Omega}=0, \ t>0,
\end{aligned}
\right.
\end{equation}
where $\textbf{v}$ is the velocity. It is well-defined as $u$ goes to zero, by the assumption that $f(0)=0$ and $f'(0)\neq 0$.

%{\color{blue}{ and $f'(0) \neq 0$. Actually the second equation of \eqref{equ:GFP-1} can be also written as
%$$ \frac{u}{f(u)}\textbf{v} =-  \partial_x \big[H'(u)+V(x)+W*u\big],$$
%where $\frac{u}{f(u)}$  makes sense, since $f'(0)\ne 0$ and $\frac{f(u)}{u}\leq C$, i.e., $\lim\limits_{u\rightarrow 0}\frac{u}{f(u)}\ne 0$.}}
 \begin{lem}\label{lem-dissip}
 If  $u(x,t)$ is the solution of \eqref{equ:GFP}, then
  $u$ satisfies the corresponding energy dissipation law
  \begin{equation}\label{equ:energylaw}
      \frac{d}{dt}E^{total}=-\Delta,
  \end{equation}
  where the total energy $$E^{total}:=\int_{\Omega}[H(u)+uV(x)]dx+\frac{1}{2}\int_{\Omega}\int_{\Omega}W(x-y)u(x)u(y)dydx,$$ and
  the entropy production $$\Delta=\int_{\Omega}\frac{u^2}{f(u)}|{\bf v}|^2dx,$$
with the velocity ${\bf v}=-\frac{f(u)}{u}\partial_x[H'(u)+V(x)+W*u]$. If $u$ satisfies the corresponding energy dissipation law \eqref{equ:energylaw}
 and a zero-flux boundary condition, then it can be shown by the Energetic Variational Approach that $u(x,t)$ solves \eqref{equ:GFP}.
\end{lem}
\noindent\textbf{Proof}:
We first prove that the energy dissipation law (\ref{equ:energylaw}) holds if $u$ is the solution of  \eqref{equ:GFP-1}.
Multiplying  by $H'(u)+V(x)+W*u$  and integrating on both sides of the first equation in \reff{equ:GFP-1}, we have
$$\int_{\Omega}[H'(u)+V(x)+W*u]\cdot\partial_t u dx=-\int_{\Omega}\partial_x(u\textbf{v})\cdot [H'(u)+V(x)+W*u]dx .$$
By integration by parts, we have
\[
\frac{d}{dt}E^{total} =\int_{\Omega}u\textbf{v}\cdot\{\partial_x[H'(u)+V(x)+W*u]\}dx = -\int_{\Omega}\frac{u^2}{f(u)}|\textbf{v}|^2dx,
\]
where the velocity $\textbf{v}=-\frac{f(u)}{u}\partial_x[H'(u)+V(x)+W*u]$ and we have used the zero-flux boundary condition $u\textbf{v}| _{\partial\Omega}=0$; cf. the system \reff{equ:GFP-1}.
%\begin{subequations}
%% \label{eq:2}
%  \begin{align}
% & \frac{d}{dt}\left(\int_{\Omega}[H(u)+uV(x)]dx+\frac{1}{2}\int_{\Omega}\int_{\Omega}W(x-y)u(x)u(y)dydx\right) \nonumber\\
% &=\int_{\Omega}u\textbf{v}\cdot\{\partial_x[H'(u)+V(x)+W*u]\}dx \nonumber\\
%  %=&-\int_{\Omega}f(u)|\partial_x[H'(u)+V(x)+W*u]|^2dx  \nonumber\\
%  &= -\int_{\Omega}\frac{u^2}{f(u)}|\textbf{v}|^2dx, \nonumber
%  \end{align}
%$$u\textbf{v}| _{\partial\Omega}=-f(u)\partial_{x}\{[H'(u)+V(x)+W*u]\} |_{\partial\Omega}=0,$$.
%\end{subequations}

Next we shall show by the EnVarA that \reff{equ:GFP-1} can be deduced from the energy dissipation law (\ref{equ:energylaw}). By mass conservation, we have
\begin{equation*}
0=\frac{d}{dt}\int_{E_t^x} u(x,t)dx =\frac{d}{dt}\int_{E_0^X} u(x(X,t),t)\det\frac{\partial x}{\partial X}dX =\int_{E_t^x} u_t+\partial_x u\cdot \textbf{v}+ u(x,t)\cdot \partial_x\textbf{v}dx, \nonumber
  \end{equation*}
 where $\textbf{v}$ denotes the velocity, $E_t^x \subset \Omega_t^x$ is the deformed configuration of an arbitrary subdomain $E_0^X \subset \Omega_0^X$, and $\det\frac{\partial x(X,t)}{\partial X}$ is the Jacobian matrix of the map: $X\rightarrow x(X,t)$.  Thus, we have
 \begin{equation}\label{ConsEq}
 u_t+\partial_x(u\textbf{v})=0.
 \end{equation}
In the Lagrangian coordinate, mass conservation leads to
 \begin{equation}\label{equ:conservationL}
 u(x(X,t),t)=\frac{u_{0}(X)}{\det\frac{\partial x(X,t)}{\partial X}},
\end{equation}
where $u_{0}(X)$ is the initial condition.

\begin{compactenum}
  \item[$\bullet$] \textbf{Least Action Principle}.
  The action functional is defined as
  \begin{equation}\label{equ:A}
  \begin{split}
  \mathcal{A}(x):=&-\int_0^{t^*}\int_{\Omega}H\Big(\frac{u_0(X)}{\partial_X x}\Big)\partial_X x\ dXdt-\int_0^{t^*}\int_{\Omega}u_0(X)V(x)dXdt\\
  &-\frac{1}{2}\int_0^{t^*}\int_{\Omega}\int_{\Omega}u_0(X)u_0(Y)W(x-y)dXdYdt,
  \end{split}
  \end{equation}
where $t^*$ is a positive number denoting the time period under consideration. %and the $y:=y(Y,t)$, $Y\in\Omega$ in Lagrangian coordinate can be  pushed by  $y$  in  Eulerian coordinate. \zhou{this y should be different from previous literature}

Based on the Least Action Principle, we have the conservative force in the Eulerian coordinate by taking the variational of $\mathcal{A}(x,t)$ with respect to $x$:
  $$ F_{con}:= \frac{\delta\mathcal{A}}{\delta x} =-u\partial_x [H'(u)+V(x)+W*u].$$
 In the Lagrangian coordinate, we have
  \begin{equation}\label{equ:ForceCon}
  F_{con} =-\partial_X\left[ \frac{u_0(X)}{\partial_X x}\cdot H'-H\right]-u_0(X)V'(x)-u_0(X)\mathcal{S}(x),\nonumber%\\
  \end{equation}
where
\begin{equation}\label{eqtraSS}
\mathcal{S}(x):=\int_{\Omega}W'(x(X,t)-y(Y,t))u_0(Y)dY.
 \end{equation}
% $$\mathcal{S}(x):=\int_{\Omega}W'(x-y)u_0(Y)dY.$$ Note that $W$ can be piecewise smooth.

%If  $W$ is  piecewise smooth at finite nodes $\bar{x}_0,\cdots, \bar{x}_{k}$, $W'(\bar{x}_i-y)$ has two values $W'_+(\bar{x}_i-y)$ and $W'_-(\bar{x}_i-y)$
%%with $\forall y:=y(Y,t)\in\Omega$,
%for $i=0,\cdots,k$.
 %if $W$ is smooth, and $$\mathcal{S}(x):=\frac{u_0}{\partial_X x}\partial_X\int_{\Omega}W(x-y)u_0(Y)dY,$$ if $W$ is integrable and non-smooth.

\item [$\bullet$] {\bf  Maximum Dissipation Law}.
     By the Maximum Dissipation Law, i.e., the Onsager's Principle, we obtain the dissipation force by taking the variation of $\frac{1}{2} \Delta$ with respect to the velocity $\textbf{v}$:
       $$F_{dis}:=\frac{\delta(\frac{1}{2}\Delta)}{\delta\textbf{v}}.$$
       Here the factor $\frac{1}{2}$ is included according to the convention that the energy dissipation $\Delta$  is always a quadratic function of certain rates, such as the velocity in the linear response theory \cite{J.W.Strutt(1873)}.

   We obtain the dissipation force

  \begin{equation}
  F_{dis}=\frac{u^2}{f(u)}{\bf v}
  \end{equation}
 and
  \begin{eqnarray}\label{equ:ForceDis}
  &&
F_{dis}=\frac{u_0^2(X)}{\partial_X x}\cdot\frac{1}{f\Big(\frac{u_0(X)}{\partial_X x}\Big)}\cdot x_t \nonumber%\\
%&&
%F_{dis}=\frac{\delta\frac{1}{2}\Delta}{\delta x_{t}}=\frac{f_{0}(X)}{x(1-x)}x_{t}.\nonumber
  \end{eqnarray}
in the Eulerian coordinate and Lagrangian coordinate, respectively.

\item [$\bullet$] {\bf  Force Balance.} By the Newton's force balance law
       $F_{con}=F_{dis},$ we have the trajectory equation
\begin{equation} \label{eqtra}
\frac{u_0^2(X)}{\partial_X x}\cdot\frac{1}{f\Big(\frac{u_0(X)}{\partial_X x}\Big)}\cdot x_t=-\partial_X \left[\frac{u_0(X)}{\partial_X x}H'\Big(\frac{u_0(X)}{\partial_X x}\Big)-H\Big(\frac{u_0(X)}{\partial_X x}\Big)\right]-u_0(X)V'(x)-u_0(X)\mathcal{S}(x),
\end{equation}
%where $\partial_X\left(G\Big(\frac{u_0}{\partial_X x}\Big)\right):=\frac{u_0}{\partial_X x}\partial_X H'$.
%where
in the Lagrangian coordinate, where $\mathcal{S}(x)$ is given in \eqref{eqtraSS}.
%\begin{equation}\label{eqtraSS}
%\mathcal{S}(x):=\int_{\Omega}W'(x(X,t)-y(Y,t))u_0(Y)dY,\ Y\in\mathcal{Q}.
% \end{equation}
In the Eulerian coordinate, we have
\begin{equation}\label{ne}
    \begin{aligned}
   \frac{u^2}{f(u)}{\bf v} =-u\partial_x \left[H'+V(x)+W*u\right].
    \end{aligned}
\end{equation}
Thus, we have by \reff{ne} that the velocity $\textbf{v}=-\frac{f(u)}{u}\partial_x[H'(u)+V(x)+W*u]$.
\end{compactenum}
Combination of \reff{ConsEq}, \reff{ne}, and the zero-flux boundary condition completes the proof. $\hfill\Box$
 %if $W$ is bound and  smooth;
 %\begin{equation}\label{eqtraSN}
 %\mathcal{S}(x):=\frac{\partial_X\left[\int_{\Omega}W(x-y)u_0(Y)dY\right]}{\partial_X x},
  %\end{equation} if $W$ is non-smooth, bounded and integral.

The trajectory $x(X,t)$ is obtained by solving \eqref{eqtra} with the initial condition
\begin{equation}\label{eqtraini}  x(X,0)=X, \ X\in\Omega, \end{equation}
and the boundary condition
\begin{equation}\label{eqtraboun} x|_{\partial\Omega}= X|_{\partial\Omega},\ t>0. \end{equation}
With the flow map $x(X,t)$, we obtain the solution $u(x,t)$ to the problem \eqref{equ:GFP} by the equation \eqref{equ:conservationL}.

\begin{rem}
Notice that the initial condition $u_0(X)>0$ on $\Omega$ is considered in Lemma \ref{lem-dissip}. The case with compactly supported initial conditions  will be treated as free boundary problems in Section \ref{sec:3}.
\end{rem}

\section{Numerical Methods for Trajectory Equation}
\label{sec:3}
In this section, we propose a structure-preserving finite difference scheme for the trajectory equation~\eqref{eqtra}.
%\begin{enumerate}[(1)]
\subsection{Time Discretization}
%\end{enumerate}
We develop a time discretization scheme using the convex splitting strategy. It follows from the convexity of $H(u)$ that the term $H\Big(\frac{u_0}{\partial_X x}\Big)\partial_X x$ is convex as well, by the assumption that $\partial_X x >0$.  The functions $V(x)$ and $W(x)$ in~\eqref{eqEne} can be split into convex part and concave part, i.e.,
$$V(x):=V_c(x){-}V_e(x),$$  $$W(x):=W_c(x){-}W_e(x),$$
where $V_c$, $W_c$, $V_e$, and $W_e$ are convex functions.
Then the trajectory equation~\eqref{eqtra} can be viewed as a gradient flow associated with the total energy
\begin{equation} \label{eqtotEn}
\begin{aligned}
E^{total}&=\int_{\Omega}H\left(\frac{u_0(X)}{\partial_X x}\right)\partial_X x+u_0(X)V(x)dX+\frac{1}{2}\int_{\Omega}\int_{\Omega}W(x-y)u_0(X)u_0(Y)dXdY\\
&:=E^{total}_c-E^{total}_e,
\end{aligned}
\end{equation}
where $$E^{total}_c:=\int_{\Omega}H\left(\frac{u_0(X)}{\partial_X x}\right)\partial_X x+u_0(X)V_c(x)dX+\frac{1}{2}\int_{\Omega}\int_{\Omega}W_c(x-y)u_0(X)u_0(Y)dXdY,$$
and $$E^{total}_e:=\int_{\Omega}u_0(X)V_e(x)dX{+}\frac{1}{2}\int_{\Omega}\int_{\Omega}W_e(x-y)u_0(X)u_0(Y)dXdY.$$
Notice that both $E^{total}_e$ and $E^{total}_c$ are convex functionals with respect to the trajectory $x(X)$.
%\begin{rem}
%Considering  the split of $W=W_c{-}W_e$,  we have the split of $$\frac{1}{2}\int_{\Omega}\int_{\Omega}W(x-y)u_0(X)u_0(Y)dYdX=\frac{1}{2}\int_{\Omega}\int_{\Omega}W_c(x-y)u_0(X)u_0(Y)dXdY
%-\frac{1}{2}\int_{\Omega}\int_{\Omega}W_e(x-y)u_0(X)u_0(Y)dXdY.$$
%\end{rem}

Based on a convex splitting technique \cite{D.J.Eyre(1998),C.H.Duan(2017)},  we propose a \emph{semi-discrete scheme} for \eqref{equ:GFP}:
% \begin{equation}\label{equ:numTime}
\[
 \begin{split}
 \frac{u^2_0(X)}{\partial_X x^n}\cdot\frac{1}{f\Big(\frac{u_0(X)}{\partial_X x^n}\Big)}\cdot\frac{x^{n+1}-x^n}{\tau}
 &=-\partial_X \left[\frac{u_0(X)}{\partial_X x^{n+1}}\cdot H'\Big(\frac{u_0(X)}{\partial_X x^{n+1}}\Big)-H\Big(\frac{u_0(X)}{\partial_X x^{n+1}}\Big)\right]\\
 &\ \ \ -u_0(X)V_c'(x^{n+1}){-}u_0(X)V_e'(x^{n})-u_0(X)\mathcal{S}_c^{n+1}{+}u_0(X)\mathcal{S}_e^{n},
 \end{split}
 \]
%\end{equation}
where $\tau:=\frac{T}{N}$, $N\in\mathbb{N}^+$,  is the time step size with the final time $T$, and
%If $W$ is bound and  smooth, then
$$\mathcal{S}_c^{n+1}:=\int_{\Omega}W_c'(x^{n+1}-y^{n+1})u_0(Y)dY $$ and $$\mathcal{S}_e^{n}:=\int_{\Omega}W_e'(x^{n}-y^{n})u_0(Y)dY.$$
 %If $W$ is non-smooth, bounded and integral, then $$\mathcal{S}_c^{n+1}:=\frac{\partial_X\left[\int_{\Omega}W_c(x^{n+1}-y^{n+1})u_0(Y)dY\right]}{\partial_X x^{n+1}},$$ and $$\mathcal{S}_e^{n}:=\frac{\partial_X\left[\int_{\Omega}W_e(x^{n}-y^{n})u_0(Y)dY\right]}{\partial_X x^{n}}.$$
%\begin{enumerate}[(2)]
\subsection{ Fully Discrete Scheme with a Positive Initial State}
%\end{enumerate}
Let $X_0$ be the left endpoint of $\Omega$ and $h=\frac{|\Omega|}{M}$  be the mesh step with $M\in\mathbb{N}^{+}$. Denote by $X_{r}=X(r)=X_0+ r h$, where $r$ takes integer or half integer values.  Let $\mathcal{E}_{M}$ and $\mathcal{C}_{M}$ be the spaces of grid functions whose domains are $\{X_{i}\ |\ i=0,...,M\}$ and $\{X_{i-\frac{1}{2}}\ |\ i=1,...,M\}$, respectively. In componentwise, these functions are identified via
$l_{i}=l(X_{i})$, $i=0,...,M$, for $l\in\mathcal{E}_{M}$, and $\phi_{i-\frac{1}{2}}=\phi(X_{i-\frac{1}{2}})$,  $i=1,...,M$, for $\phi\in\mathcal{C}_{M}$. Without ambiguity,  we denote by $X=\{X_{i}\ |\ i=0,...,M\}$ for $X\in \mathcal{E}_{M}$.

We define difference operators $D_{h}: \mathcal{E}_{M}\rightarrow\mathcal{C}_{M}$, $d_{h}: \mathcal{C}_{M}\rightarrow\mathcal{E}_{M}$, and $\widetilde{D}_h: \mathcal{E}_{M}\rightarrow\mathcal{E}_{M}$ by
\begin{align}\label{equ:dif1}
& (D_{h}l)_{i-\frac{1}{2}}= (l_{i}-l_{i-1})/h,\ i=1,...,M, \\
&  (d_{h}\phi)_{i}= (\phi_{i+\frac{1}{2}}-\phi_{i-\frac{1}{2}})/h,\ i=1,...,M-1,\\
&(\widetilde{D}_hl)_{i}=\left\{
\begin{array}{lcl}
(l_{i+1}-l_{i-1})/2h, &\mbox{$ i=1,...,M-1$},\\
(l_{i+1}-l_{i})/h, &\mbox{$ i=0$},\\
(l_{i}-l_{i-1})/h, &\mbox{$ i=M$}.
% (4l_{i+1}-l_{i+2}-3l_{i})/2h,&  i=0,  \\
% (l_{i-2}-4l_{i-1}+3l_{i})/2h, &  i=M.
\end{array}\right.
\end{align}
Let $l$, $g\in\mathcal{E}_M$ and $\phi$, $\varphi\in\mathcal{C}_M$. We define the \emph{inner product} on space $\mathcal{E}_M$ and $\mathcal{C}_M$ by
\begin{eqnarray}
  %\|l \|_2 =\sqrt{ \left\langle l , l \right\rangle } ,  \quad \mbox{with} \quad
  &&
 \left\langle l , g \right\rangle_{\mathcal{E}}  := h \left( \frac12 l_0 g_0 + \sum_{i=1}^{M-1}   l_{i} g_{i}
  + \frac12 l_M g_M \right),\\
  \label{FD-inner product-1}
%\end{eqnarray}
%Similarly, for any grid functions $\phi$ and $\varphi$, evaluated at the stagger grid point $X_{i+1/2}$, the discrete $L^2$ inner product and norm becomes
%\begin{eqnarray}
&&
 % \| \phi\|_2 = \sqrt{ \left\langle\phi ,  \phi \right\rangle_e } ,  \quad \mbox{with} \quad
  \left\langle \phi , \varphi \right\rangle_{\mathcal{C}}
  := h \sum_{i=0}^{M-1}   \phi_{i+\frac 12} \varphi_{i+\frac 12} .
  \label{FD-inner product-2}
\end{eqnarray}
It is easy to verify the following summation by parts formula:
\begin{equation}\label{FD-inner product-3-1}
   \left\langle l , d_h \phi \right\rangle_{\mathcal{E}}
  = -\left\langle D_h l , \phi \right\rangle_{\mathcal{C}},
  \mbox{\ with $l_0 = l_M =0$, $\phi\in\mathcal{C}_{M}$, $l\in\mathcal{E}_M$}.
\end{equation}
%When $l(X)$ is discontinuous at  finite many breakpoints $ \bar{z}_0,\cdots, \bar{z}_{\bar{K}}$, i.e., $l$ has two values at each point,  the  \emph{inner product}  can be defined on $\mathcal{K}:=\{\bar{z}_0,\cdots, \bar{z}_{\bar{K}}\}$ as
%\begin{eqnarray}  \label{FD-inner product-3}
%&&
% \left\langle l, g \right\rangle_{\mathcal{K}}=\sum\limits_{i=0}^{\bar{K}} \left\langle l, g \right\rangle_{i},
%  {\mbox{with}\ } \left\langle l, g \right\rangle_{i}=\frac{h}{2} \left( l^-_ i g_i +  l^+_{i} g_{i} \right).%\nonumber
%\end{eqnarray}
%
%\begin{rem}
%By the definition \eqref{FD-inner product-3}, the precision  of  integral can be held. \zhou{I do not understand}
%\end{rem}

Let $\mathcal{Q}:=\{l \in\mathcal{E}_{M}\ |\ l_{i-1}<l_{i},\ 1\leq i\leq M;\ l_{0}=X_0,\ l_{M}=X_M\}$ be an admissible set, in which particles are arranged in the order without twisting or exchanging. Its boundary set is given by $\partial\mathcal{Q}:=\{l \in\mathcal{E}_{M}\ |\ l_{i-1}\leq l_i,\ 1\leq i\leq M,\ \text{and exists } i_0\in\{1,\cdots,M\} \text{~such that} \ l_{i_0}=l_{i_0-1};\ l_{0}=X_0,\ l_{M}=X_M\}$. Clearly, $\bar{\mathcal{Q}}:=\mathcal{Q}\cup\partial\mathcal{Q}$ is a closed convex set.

The {\bf  fully discrete scheme} is formulated as follows: Given $x^{n} \in\mathcal{Q}$, find $x^{n+1}=(x^{n+1}_{0},...,x^{n+1}_{M})\in\mathcal{Q}$ such that
\begin{equation}\label{equ:numnum}
\begin{split}
\frac{u^2_{0_i}}{\widetilde{D}_h x^n_i}\cdot\frac{1}{f\Big(\frac{u_{0_i}}{\widetilde{D}_h x^n_i}\Big)}\cdot\frac{x^{n+1}_i-x^n_i}{\tau}
 &=-d_h \left[\frac{u_0}{D_h x^{n+1}}\cdot H'\Big(\frac{u_0}{D_h x^{n+1}}\Big)-H\Big(\frac{u_0}{D_h x^{n+1}}\Big)\right]_i\\
 &\ \ \ \ -u_{0_i} V_c'(x^{n+1}_i){+}u_{0_i} V_e'(x^{n}_i)-u_{0_i}\hat{\mathcal{S}}_{c_{i}}^{n+1}{+}u_{0_i} \hat{\mathcal{S}}_{e_i}^n,\ 1\leq i\leq M-1,
 \end{split}
\end{equation}
with boundary conditions
\begin{equation}\label{equ:numbou}
 x_0^n=X_0, \ x_M^n=X_M.
\end{equation}
Here $\hat{\mathcal{S}}_{c_{i}}^{n+1}$ and $\hat{\mathcal{S}}_{e_i}^n$ are given by
%\begin{itemize}
%\item {\color{blue}If $W$ is smooth and integral,}

 \begin{equation}\label{eqnumSSconv}
    %\mathcal{S}_{c_{i}}^{n+1}:=\sum\limits_{j=1}^{M-1}\frac{h}{6}\big[W'_c(x^{n+1}_i-y^{n+1}_{j-1})u_{0_{j-1}}+4W'_c(x^{n+1}_i-y^{n+1}_{j})u_{0_j}
    %+W'_c(x^{n+1}_i-y^{n+1}_{j+1})u_{0_{j+1}}\big],
\hat{\mathcal{S}}_{c_{i}}^{n+1}:=\left\langle W'_c(x^{n+1}_i-y^{n+1}),u_0(Y)\right\rangle_{\mathcal{E}},
    \end{equation}
and
\begin{equation}\label{eqnumSSconc}%\mathcal{S}_{e_{i}}^{n}:=\sum\limits_{j=1}^{M-1}\frac{h}{6}\big[W'_e(x^{n}_i-y^{n}_{j-1})u_{0_{j-1}}+4W'_e(x^{n}_i-y^{n}_{j})u_{0_j}
    %+W'_e(x^{n}_i-y^{n}_{j+1})u_{0_{j+1}}\big].
\hat{\mathcal{S}}_{e_{i}}^{n}:=\left\langle W'_e(x^{n}_i-y^{n}),u_0(Y)\right\rangle_{\mathcal{E}},
    \end{equation}
with $x^{n+1}_i:=x(X_i,t^{n+1})$ and $y^{n+1}:=y(Y,t^{n+1})$.

%\item If  $W$ is non-smooth and integral,
%\begin{equation}\label{eqnumSNconv}
%\mathcal{S}_{N_{c_{i}}}:=\frac{(\bar{S}^{n+1}_{c_{i+1}}-\bar{S}^{n+1}_{c_{i-1}})/2h}
%{(x^{n+1}_{i+1}-x^{n+1}_{i-1})/2h},\end{equation}
%with
%$$\bar{S}^{n+1}_{c_{i}}:=\sum\limits_{j=1}^{M-1}\frac{h}{6}\big[W_c(x^{n+1}_i-y^{n+1}_{j-1})u_{0_{j-1}}+4W_c(x^{n+1}_i-y^{n+1}_{j})u_{0_j}
%    +W_c(x^{n+1}_i-y^{n+1}_{j+1})u_{0_{j+1}}\big],$$
%
%\begin{equation}\label{eqnumSNconc}
%\mathcal{S}_{N_{e_i}}^{n}:=\frac{(\bar{S}^{n}_{e_{i+1}}-\bar{S}^{n}_{e_{i-1}})/2h}
%{(x^{n}_{i+1}-x^{n}_{i-1})/2h},\end{equation}
%with
%$$\bar{S}^{n}_{e_{i}}:=\sum\limits_{j=1}^{M-1}\frac{h}{6}\big[W_e(x^{n}_i-y^{n}_{j-1})u_{0_{j-1}}+4W_e(x^{n}_i-y^{n}_{j})u_{0_j}
%    +W_e(x^{n}_i-y^{n}_{j+1})u_{0_{j+1}}\big].$$
%\end{itemize}
We develop a Newton's iteration method to solve the nonlinear difference equations~\eqref{equ:numnum}. Define the following convex functional
\begin{equation}\label{minEnergy}
\begin{split}
 J(z):=&\frac{1}{2\tau}\left\langle\frac{u_0^2}{\widetilde{D}_h x^n}\cdot\frac{1}{f\Big(\frac{u_0}{\widetilde{D}_h x^n}\Big)}\cdot(z-x^n),z-x^n\right\rangle_{\mathcal{E}}+\left\langle H\Big(\frac{u_0}{D_h z}\Big)D_h z,\textbf{e}\right\rangle_{\mathcal{C}}\\
&+\left\langle u_0(X),V_c(z)\right\rangle_{\mathcal{E}}+\frac{1}{2}\left\langle\Big\langle W_c(z-y),u_0(Y)\Big\rangle,u_0(X)\right\rangle_{\mathcal{E}}\\
&{-}\left\langle u_0(X)V'_e(x^n),z\right\rangle_{\mathcal{E}}{-}\left\langle u_0(X)\Big\langle W'_e(x^n-y^n),u_0(Y)\Big\rangle,z\right\rangle_{\mathcal{E}},
\end{split}
\end{equation}
where  $x^{n},y^n\in\mathcal{Q}$ are coordinates of particles at time $t^{n}$, $n=0,\cdots,N-1$, and $\textbf{e}$ is a vector with each element being one.

\noindent{\bf Newton's iteration}. \ \ Set $x^{n+1, 0}= x^n$. For $k=0,1,2,\cdots$, update $x^{n+1,k+1}=x^{n+1,k} +\delta_x$, where $\delta_x$ solves equations
\begin{equation}\label{equ:numnumNew}
\begin{split}
&\frac{u^2_{0_i}}{\widetilde{D}_h x^n_i}\cdot\frac{1}{f\Big(\frac{u_{0_i}}{\widetilde{D}_h x^n_i}\Big)}\cdot\frac{x^{n+1,k}_i+\delta_{x_i}-x^n_i}{\tau}\\
 &=-d_h \left[\frac{u_0}{D_h x^{n+1,k}}\cdot H'\Big(\frac{u_0}{D_h x^{n+1,k}}\Big)-H\Big(\frac{u_0}{D_h x^{n+1,k}}\Big)\right]_i\\
 &\ \ \ +d_h \left[\Big(\frac{u_0^2}{(D_h x^{n+1,k})^3} H''\Big)D_h \delta_x\right]_i-u_{0_i}V_c'(x^{n+1,k}_i)-u_{0_i} V_c''(x^{n+1,k}_i)\delta_{x_i}\\
 &\ \ \ {+}u_{0_i} V_e'(x^{n}_i)-u_{0_i}\hat{\mathcal{S}}_{c_{i}}^{n+1,k}-u_{0_i}(\hat{\mathcal{S}}'_c)_i^{n+1,k}\delta_{x_i}{+}u_{0_i}\hat{
 \mathcal{S}}_{e_i}^n,\ 1\leq i\leq M-1,
 \end{split}
\end{equation}
with boundary conditions $\delta_{x_0}=\delta_{x_M}=0$, where
$$(\hat{\mathcal{S}_c}')_{i}^{n+1,k}:=\left\langle W''_c(x^{n+1,k}_i-y^{n+1,k}),u_0(Y)\right\rangle_{\mathcal{E}}.$$
\begin{rem}
To guarantee the convergence of the Newton's iterations on $\mathcal{Q}$ theoretically, the damped Newton's iteration \cite{Y.Nesterov(1994)} can be used \cite{C.H.Duan(2017),C.H.Duan(2018)}. However, the proposed Newton's iteration method is more efficient and converges robustly in many examples.
\end{rem}

After solving~\eqref{equ:numnumNew}, we finally get numerical density $\{u^{n}_i\}_{i=0}^M$, $n=1,\cdots,N$, from \eqref{equ:conservationL} by
\begin{align}
& u_{i}^{n}= \frac{ u_{0}(X_{i}) } {(x_{i+1}^{n}-x_{i-1}^{n})/(2h)},\ 1\leq i\leq M-1, \ \mbox{and} \label{numdist}\\
& u_{0}^{n}=\frac{u_{0}(X_{0})} {(x_{1}^{n}-x_{0}^{n})/h}, \ u_{M}^{n}=\frac {u_{0}(X_{M})}{(x_{M}^{n}-x_{M-1}^{n})/h}. \label{numdist2}
\end{align}
%{\color{blue}{\begin{rem}
%
% In the discrete scheme \eqref{equ:numnum}, we find the solution $x^{n+1}\in\mathcal{Q}$, i.e., $x_{i}<x_{i+1}$, $i=0,\cdots,M-1$, so that the particles never touch each other.  Hence \eqref{numdist} - \eqref{numdist2} makes sense.
%
% In  some cases, however,  the distance between a certain particle $x^{n+1}_{i_0}$ and its neighbor $x^{n+1}_{i_0+k}$ or/and $x^{n+1}_{i_0-k}$, $k=1,2,\cdots$,  may be too close to distinguish from each other under the machine precision, they are bundled up and will be regarded as one particle, which is the signal that the numerical Dirac delta happens. The density at $x^{n+1}_{i_0}$ should be changed and can be approximated at the machine precision. The more details can be found in {\bf{ Example 4}} in Sec. \ref{sec:5} or in \cite{C.H.Duan(2017)}.
%\end{rem}}}
\begin{thm}\label{thm:Con+Pos}
The numerical solution $\{u^{n}_i\}_{i=0}^M$, $n=1, \cdots, N$, obtained from \eqref{numdist}-\eqref{numdist2}, ensures mass  conservation and positivity of solution with an initial condition $u_0> 0$.
\end{thm}
\noindent\textbf{Proof:}  Let the initial mass carried by each particle $x_i^0=X_i$ be
\begin{equation}\label{init-mass}
m_i^0= h u_0(X_i), \ 1<i<M; \ m_0^0 = \frac h2 u_0(X_0); ~\text{	and} \ m_M^0 = \frac h2 u_0(X_M).
\end{equation}
Define the mass carried by each particle $x_i^n$ as
\begin{equation}\label{mass}
m_i^n= \frac{x_{i+1}^n-x_{i-1}^n} 2  u_i^n, \ 1<i<M; \ \ m_0^n = \frac {x_1^n-x_0^n}2 u_0^n; ~\text{	and}\ m_M^n = \frac {x_M^n-x_{M-1}^n}2 u_M^n.
\end{equation}
It follows from \eqref{numdist}-\eqref{numdist2} that
$$m_i^n \equiv m_i^0, \ \ 0\leq i\leq M, \ n=1,2,\cdots N.$$
Since $x^{n} \in \mathcal{Q}$, we know that
$u_i^n > 0, \ \ 0\leq i\leq M, \ n=1,2,\cdots N$. $\hfill\Box$
\subsection{Numerical Scheme  for Problems with Free Boundaries}
For the initial data with a compact support in $\Omega$, we define the left and right interfaces as
\begin{equation} \xi_{1}^t:=\inf\{x\in\Omega:u(x,t)>0,t\geq 0\}, \label{equ:boundaryleft} \end{equation}
\begin{equation} \xi_{2}^t:=\sup\{x\in\Omega:u(x,t)>0,t\geq 0\}. \label{equ:boundaryright} \end{equation}
 Let $\Gamma^{t}:=[\xi_{1}^t, \xi_{2}^t]\subset\Omega$. %For this kind of problems, all the trajectories start from the initial support $\Gamma^0\subsetneqq\Omega$.
Assume $H(\cdot)\in C^2([0,+\infty))$ and $f(\cdot)\in C^1([0, +\infty))$ with $f'(0)\neq 0$. %(We use the standard notation for Sobolev spaces; cf.\\cite{Adams75}.)
% Let $\lim\limits_{u\rightarrow 0}\frac{f(u)}{u}= f'(0)$.

 We shall solve the initial-boundary value problem~\eqref{eqtra} with the boundary condition
  %\begin{eqnarray}
% && x_t=-f'\Big(\frac{u_0(X)}{\partial_X x}\Big) \frac{\partial_X u_0(X)}{(\partial_X x)^2} H''-f'\Big(\frac{u_0(X)}{\partial_X x}\Big)V'(x)-f'\Big(\frac{u_0(X)}{\partial_X x}\Big)\mathcal{S}(x),\ X\in\partial\Gamma^0, \label{eqtraFbou}\\
  % && f(x,0)=f_{0}(x),\ x\in\Omega \\
  %   && x(X,0)=X, \ X\in\Gamma^0,  \label{eqtraFini}
   % \end{eqnarray}
\begin{equation}\label{eqtraFbou}
x_t=-f'(0) \frac{\partial_X u_0(X)}{(\partial_X x)^2} H''\left(\frac{u_0(X)}{\partial_X x}\right)-f'(0)V'(x)-f'(0)\mathcal{S}(x),\ X\in\partial\Gamma^0,
\end{equation}
and the initial condition
\begin{equation}\label{eqtraFini}
x(X,0)=X, \ X\in\Gamma^0,
\end{equation}
where time $t>0$ and $\mathcal{S}(x)$ is given by \eqref{eqtraSS}.
%The first term may be written as  $$f'{\color{blue}(0)}\cdot \frac{\partial_X g( u_0(X))}{h(\partial_X x)},$$ with $H''=\frac{r(u_0(X))}{\bar{r}(\partial_X x)}$.
Here the boundary condition~\eqref{eqtraFbou} on the free boundary is obtained from~\eqref{eqtra} with the condition $u_0(X)|_{\partial\Gamma}=0$.% {\color{blue} and the condition also holds in $H''$...}%

We partition the interval $\Gamma^0$ into equal subintervals with $X_i=\xi_{1}^0+ih$, $0\leq i \leq M$, where $M$ is an integer number and the grid spacing $h:=(\xi_{2}^0-\xi_{1}^0)/M$. Given the initial state $u_0(X)$ with a compact support $\Gamma^0$, the numerical solution for the trajectory equation is obtained by solving~\eqref{equ:numnum} for $1\leq i\leq M-1$, and discrete boundary conditions
%The \textbf{fully discrete} schemes is: Given  the initial state $u_0(X)$ with a compact support $\Gamma^0$ and $\{x_i^n\}_{i=0}^M$, find $\{x_i^{n+1}\}_{i=0}^M$ such that
%\begin{equation}
% \begin{split}\label{eqnumbou}
% \frac{x^{n+1}_i-x^{n}_i}{\tau}=&-f'\Big(\frac{u_0(X)}{\widetilde{D}_h  x^n}\Big)\cdot H''\left(\frac{u_0(X)}{\widetilde{D}_h x^{n+1}_i}\right)\cdot\frac{\widetilde{D}_h u_0(X)}{(\widetilde{D}_h x^{n+1}_i)^2}-f'\Big(\frac{u_0(X)}{\widetilde{D}%_h  x^n}\Big)\cdot V'_c(x^{n+1}_i) \nonumber \\
% &{+}f'\Big(\frac{u_0(X)}{\widetilde{D}_h  x^n}\Big)V'_e(x^n_i)-f'\Big(\frac{u_0(X)}{\widetilde{D}_h  x^n}\Big)\cdot\hat{\mathcal{S}}^{n+1}_{c_i}{+}f'\Big(\frac{u_0(X)}{\widetilde{D}_h  x^n}\Big)\hat{\mathcal{S}}^{n}_{e_i},\ \ i=0,M,
 % \end{split}
  %\end{equation}
\begin{equation}
 \begin{split}\label{eqnumbou}
 \frac{x^{n+1}_i-x^{n}_i}{\tau}=&-f'(0)\cdot H''\left(\frac{u_0(X)}{\widetilde{D}_h x^{n+1}_i}\right)\cdot\frac{\widetilde{D}_h u_0(X)}{(\widetilde{D}_h x^{n+1}_i)^2}-f'(0) V'_c(x^{n+1}_i) \\
 &\quad{+}f'(0)V'_e(x^n_i)-f'(0)\hat{\mathcal{S}}^{n+1}_{c_i}{+}f'(0) \hat{\mathcal{S}}^{n}_{e_i}~ \text{ for } i=0 \text{ and }M,
  \end{split}
  \end{equation}
  where
  $\hat{\mathcal{S}}^{n+1}_{c_i}$ and $\hat{\mathcal{S}}^{n}_{e_i}$ are given by \eqref{eqnumSSconv} and \eqref{eqnumSSconc}, respectively.
  %$$\mathcal{S}^{n+1}_{N_{e_i}}:=\frac{4\bar{\mathcal{S}}^n_{e_{i+1}}-
%  \bar{\mathcal{S}}^n_{e_{i+2}}-3\bar{\mathcal{S}}^n_{e_{i+1}}}{4x^n_{i+1}-
%  x^n_{i+2}-3x^n_{i}},\ i=0,$$
%  $$\mathcal{S}^{n+1}_{N_{e_i}}:=\frac{\bar{\mathcal{S}}^n_{e_{i-2}}-4
%  \bar{\mathcal{S}}^n_{e_{i-1}}+3\bar{\mathcal{S}}^n_{e_{i}}}{x^n_{i-2}-
%  4x^n_{i-1}+3x^n_{i}},\ i=M,$$
%  and $\mathcal{S}^{n+1}_{N_{c_i}}$ has the similar formula.
The whole nonlinear system is again solved by a Newton's iteration method.

The waiting time phenomenon, i.e., a free boundary remains stationary during $(0,t^*)$ for $0<t^*<\infty$, may occur under some conditions. Such a phenomenon is common and important for the porous medium type of equations \cite{J.L.Vazquez(2007),D.G.Aronson(1983)}. However, the development of efficient algorithms for the calculation of  waiting time is challenging. Based on the Energetic Variational Approach, we have proposed an algorithm to calculate the waiting time for porous medium equations in \cite{C.H.Duan(2018)}. In this work, we extend the algorithm to consider more general problems described by the nonlinear Fokker--Planck equations that have the waiting time phenomenon.

As the right hand side of the equation \eqref{eqtraFbou} vanishes, the free boundaries remain stationary.
%That means the right term of \eqref{eqtraFbou} equals zero at first. We success the algorithm of the waiting time in \cite{C.H. Duan(2018)}.
Without loss of generality, we consider the left boundary. Analogous results can be obtained for the right boundary.  The waiting time can be characterized by
\begin{equation*}%\label{WTC_WT}
t^*:=\inf\Big\{t>0: x_t  < 0,\  \mbox{as}\  X\to \xi_1^0\Big\}.
\end{equation*}
% Let
 %\begin{equation}\label{WTC}
% \begin{split}
 %\mathcal{B}_h^n:=&f'\Big(\frac{u_0(X)}{\widetilde{D}_h  x^n_{h,0}}\Big)\cdot H''\left(\frac{u_0(X)}{\widetilde{D}_h x^{n}_{h,0}}\right)\cdot\frac{\widetilde{D}_h u_0(X)}{(\widetilde{D}_h x^{n}_{h,0})^2}+f'\Big(\frac{u_0(X)}{\widetilde{D}_h  x^n_{h,0}}\Big)\cdot V'_c(x^{n}_{h,0})\\
% &{-}f'\Big(\frac{u_0(X)}{\widetilde{D}_h  x^n_{h,0}}\Big)\cdot V'_e(x^{n}_{h,0})+f'\Big(\frac{u_0(X)}{\widetilde{D}_h  x^n_{h,0}}\Big)\cdot \hat{\mathcal{S}}^{n+1}_{c_i}{-}f'\Big(\frac{u_0(X)}{\widetilde{D}_h  x^n_{h,0}}\Big)\cdot\hat{\mathcal{S}}^{n}_{e_i},
% \end{split}\end{equation}
We define
\begin{equation*}%\label{WTC}
 \begin{split}
 \mathcal{B}_h^n:=&f'(0) H''\left(\frac{u_0(X)}{\widetilde{D}_h x^{n}_{h,0}}\right)\cdot\frac{\widetilde{D}_h u_0(X)}{(\widetilde{D}_h x^{n}_{h,0})^2}+f'(0) V'_c(x^{n}_{h,0})\\
 &\quad{-}f'(0) V'_e(x^{n}_{h,0})+f'(0)  \hat{\mathcal{S}}^{n+1}_{c_0}{-}f'(0)\hat{\mathcal{S}}^{n}_{e_0},
 \end{split}\end{equation*}
where   $x^{n}_h=(x_{h,0}^{n},\cdots,x_{h,M}^{n})$ denotes the numerical solution at time $t^{n}$, $n=0,\cdots,N$, with a grid spacing $h$.
 The waiting time $t^*_h$  is numerically determined by the criterion \cite{C.H.Duan(2018)}
 \begin{equation}\label{WTC2_1}
t^*_h:=\min\Big\{t^n: \Big|\frac{\mathcal{B}^n_{2h}}{\mathcal{B}_h^n}\Big|\leq 1 \Big\}.
\end{equation}
% The details of the algorithm in {\bf Example 2}.

\section{Analysis Results}
\label{sec:4}
In this section, we perform numerical analysis on the numerical scheme~\eqref{equ:numnum}, including unique solvability in the admissible set, total energy dissipation, and convergence rate.

%The inverse inequality holds:
%\begin{equation}\label{equ:inverse}
%\|l\|_{\infty}\leq C_m \frac{\|l\|_2}{h^{1/2}}, \ \  l\in\mathcal{E}_M,
%\end{equation}
%where
%$$\|l\|_{\infty}:=\max\limits_{0\leq i\leq M} \{l_i\}, \mbox{\ \ and\ \ \ }
%$$\|l\|^2:=\left\langle l,l \right\rangle_{\mathcal{E}},$$
%and $C_m$ is a positive constant independent of the grid spacing.

%First we prove that there exists a unique solution in admissible  set $\mathcal{Q}$.
\begin{thm}
\label{lem:unique}
%Suppose the interaction kernel $W$ is piecewise smooth.
Assume that $H(\cdot)$, $V_c(\cdot)$, $V_e(\cdot)$, $W_c(\cdot)$, and $W_e(\cdot)$ are piecewise $C^1$ functions. If the initial state $u_0(X)\in\mathcal{E}_{M}$ is positive for $X\in \mathcal{Q}$, then the numerical scheme~\eqref{equ:numnum} is uniquely solvable in  $\mathcal{Q}$. %of the trajectory equation based on the total energy $f\ln f$ is  \eqref{eqtranum1}  and
\end{thm}
\noindent\textbf{Proof:} We first consider the following optimization problem
\begin{equation*}
\min\limits_{z\in\bar{\mathcal{Q}}} J(z),
\end{equation*}
where the functional $J$ is given in~\eqref{minEnergy}. Since $J(z)$ is a convex functional on the closed convex set $\bar{\mathcal{Q}}$,  there exists a unique minimizer  $x\in\bar{\mathcal{Q}}$. For any $z\in\partial\mathcal{Q}$, there exists some $i>0$ such that $(D_h z)_{i-1/2} = (z_i-z_{i-1})/h=0$, indicating that $J(z)=+\infty$. Therefore, we have that the minimizer $x\in\mathcal{Q}$.  Next we shall prove  that $x$ is the minimizer of $J(z)$  if and only if it solves the equations~\eqref{equ:numnum}. We then claim that the fully discrete scheme~\eqref{equ:numnum}  has a unique solution.

Suppose $x\in\mathcal{Q}$ is the minimizer of $J(z)$. Since $\mathcal{Q}$ is an open convex set, for any $ z\in\bar{\mathcal{Q}}$, there exists a sufficiently small $\varrho_0>0$, such that, for $\varrho\in(-\varrho_0, \varrho_0)$, $x+\varrho (z-x) \in \mathcal{Q}$.  Then $j(\varrho) := J(x+\varrho (z-x))$ achieves its minimum at $\varrho =0$.  So we have $j'(0)=0$. By the summation by parts \reff{FD-inner product-3-1}, we obtain for any $z\in\bar{\mathcal{Q}}$ that
\begin{equation}\label{eqTraSum}
\begin{split}
    &\frac{1}{\tau}\left\langle\frac{u_0^2}{\widetilde{D}_h x^n}\cdot\frac{1}{f\Big(\frac{u_0(X)}{\widetilde{D}_h x^n}\Big)}\cdot(x-x^n),z-x\right\rangle_{\mathcal{E}}+
    \left\langle d_h\left(\frac{u_0}{D_h x}H'-H\right),z-x\right\rangle_{\mathcal{E}}\\
    &+\left\langle u_0(X)V'_c(x),z-x\right\rangle_{\mathcal{E}}+\left\langle u_0(X)\Big\langle W'_c(x-y),u_0(Y)\Big\rangle,z-x\right\rangle_{\mathcal{E}}\\
    &{-}\left\langle u_0(X)V'_e(x^n),z-x\right\rangle_{\mathcal{E}}{-}\left\langle u_0(X)\Big\langle W'_e(x^n-y^n),u_0(Y)\Big\rangle,z-x\right\rangle_{\mathcal{E}}=0.
    \end{split}
\end{equation}
This implies that $x$ solves the equations~\eqref{equ:numnum}.

Conversely, let $x$ be a solution to the fully discrete scheme~\eqref{equ:numnum}. We shall prove that $x$ is the minimizer of $J(z)$ on $\bar{\mathcal{Q}}$. For any $ z\in\partial\mathcal{Q}$, we have  $J(z)\geq J(x)$ due to $J(z)=+\infty$. For  any $z\in\mathcal{Q}$, taking an inner product of \eqref{equ:numnum} with $z-x$ and using summation by parts lead to \eqref{eqTraSum}. By the convexity of $J(z)$ and \eqref{eqTraSum}, we have for any $z\in\mathcal{Q}$ that
  \begin{align}
  J(z)=J(x+(z-x)) &\geq J(x).%+\frac{1}{\tau}\left\langle\frac{u_0^2}{\widetilde{D}_h x^n}\cdot\frac{1}{f\Big(\frac{u_0(X)}{\widetilde{D}_h x^n}\Big)}\cdot(x-x^n),z-x\right\rangle_{\mathcal{E}}
%     \nonumber\\
   %  &\ \ \ \ \ \ \ \ \ \ +\frac{1}{2\tau}\left\langle\frac{u_0^2}{\widetilde{D}_h x^n}\cdot\frac{1}{f\Big(\frac{u_0(X)}{\widetilde{D}_h x^n}\Big)}\cdot(z-x),z-x\right\rangle_{\mathcal{E}} \nonumber  \\
  %   &\ \ \ \ \ \ \ \ \ \ +\Big\langle-\frac{u_0}{D_h x}H'+H,D_h(z-x) \Big\rangle_{\mathcal{C}} \nonumber  \\
  %   &\ \ \ \ \ \ \ \ \ \ +\left\langle u_0(X)(V'_c(x){-}V'_e(x^n)),z-x\right\rangle_{\mathcal{E}}\nonumber  \\
   %  &\ \ \ \ \ \ \ \ \ \ + \left\langle u_0(X)\Big\langle W'_c(x-y){-}W'_e(x^n-y^n),u_0(Y)\Big\rangle,z-x\right\rangle_{\mathcal{E}}\nonumber  \\
  %&\geq J(x).
  \label{equ:vp3}
  \end{align}
This completes the proof. $\hfill\Box$

We next consider discrete energy dissipation of the numerical scheme~\eqref{equ:numnum}. The discrete total energy $E_{N}:\mathcal{Q}\rightarrow \mathbb{R}$ corresponding to \eqref{eqtotEn} is defined by
\begin{equation}\label{equ:totalenergyNum}
E_{N}(x):=\left\langle H\left(\frac{u_0(X)}{D_h x}\right),D_h x\right\rangle_{\mathcal{C}}+\left\langle u_0(X),V(x)\right\rangle_{\mathcal{E}}+\frac{1}{2}\left\langle u_0(X), \left\langle W(x-y), u_0(Y)\right\rangle\right\rangle_{\mathcal{E}}.
\end{equation}
Analogously, we have the splitting $E_{N}(x)=E_{N,c}(x)-E_{N,e}(x)$, where
$$ E_{N,c}(x):=\left\langle H\left(\frac{u_0(X)}{D_h x}\right),D_h x\right\rangle_{\mathcal{C}}+\left\langle u_0(X),V_c(x)\right\rangle_{\mathcal{E}}+\frac{1}{2}\left\langle u_0(X), \left\langle W_c(x-y), u_0(Y)\right\rangle\right\rangle_{\mathcal{E}},$$
$$E_{N,e}(x):=\left\langle u_0(X),V_e(x)\right\rangle_{\mathcal{E}}{+}\frac{1}{2}\left\langle u_0(X), \left\langle W_e(x-y), u_0(Y)\right\rangle\right\rangle_{\mathcal{E}}.$$
Here both $ E_{N,c}(x)$ and $E_{N,e}(x)$ are convex and their first variations are given by
\begin{equation} \label{variation_convex} \delta_{x} E_{N,c}(x^{n+1})=d_h \left[\frac{u_0(X)}{D_h x^{n+1}}\cdot H'\Big(\frac{u_0(X)}{D_h x^{n+1}}\Big)-H\Big(\frac{u_0(X)}{D_h x^{n+1}}\Big)\right]+u_0(X)V'_c(x^{n+1})+\mathcal{S}_{N_c}^{n+1},\end{equation}
 \begin{equation}\label{variation_concave}
 \delta_{x}E_{N,e}(x^n)=u_0(X)V'_e(x^{n}){+}\mathcal{S}_{N_e}^n.\end{equation}
%N ow we show that there exists a convex splitting for the discrete total energy $E_{N}$.
\begin{thm}\label{thm:ene}
Suppose $x^{n}=(x_{0}^{n},...,x^{n}_{M}) \in\mathcal{Q}$ be the solution to the scheme \eqref{equ:numnum} at time $t^{n}$. Denote by $E_N^n:=E_{N}(x^n).$ Then the discrete  energy dissipation law holds, i.e., \begin{equation}\label{eqi:ene_total}
\frac{E_{N}^{n+1}-E_{N}^n}{\tau } \leq - \left\langle\frac{u_0^2(X)}{\widetilde{D}_h x^n}\cdot\frac{1}{f\Big(\frac{u_0(X)}{\widetilde{D}_h x^n}\Big)}\cdot\frac{x^{n+1}-x^{n}}
{\tau},\frac{x^{n+1}-x^{n}}{\tau}\right\rangle_{\mathcal{E}} \leq 0, \ \ n=0,1,\cdots.\end{equation}
This is the discrete counterpart of the dissipation law \eqref{equ:energylaw}.% in% Lemma \ref{lem-dissip}.
\end{thm}
It is easy to prove the theorem by using the convexity of $E_{N,c}^{n}$ and $E_{N,e}^{n}$. Similar ideas can be found in the works \cite{C.H.Duan(2017),C.H.Duan(2018),C.H.Duan(2019)}.

%\noindent$\textbf{Proof.}$
%Thanks to the convexity of $E_{N,c}^{n}$ and $E_{N,e}^{n}$, we have
%\begin{align} & E_{N,c}(x^{n})-E_{N,c}(x^{n+1})\geq
%  \left\langle\delta_{x} E_{N,c}(x^{n+1}),x^{n}-x^{n+1}\right\rangle_{\mathcal{E}},\notag \\
%& E_{N, e}(x^{n+1})-E_{N,e}(x^{n})\geq\left\langle\delta_{x} E_{N,e}(x^n),x^{n+1}-x^{n}\right\rangle_{\mathcal{E}}. \notag
%\end{align}
%It follows from \eqref{variation_convex}, \eqref{variation_concave}, and \eqref{equ:numnum} that
%\begin{subequations}
%\begin{align}
%E_{N}(x^{n+1})-E_{N}(x^n)&=(E_{N,c}(x^{n+1})-E_{N,e}(x^{n+1}))-(E_{N,c}(x^{n})-E_{N,e}(x^{n}))\nonumber\\
 %&=(E_{N,c}^{n+1}-E_{N,c}^{n})-(E_{N,e}^{n+1}-E_{N,e}^{n})\nonumber \\
% &\leq\left\langle\delta_{x}E_{N,c}(x^{n+1})-\delta_{x} E_{N,e}(x^n),x^{n+1}-x^{n}\right\rangle_{\mathcal{E}} \nonumber\\
%&= - \left\langle\frac{u_0^2(X)}{\widetilde{D}_h x^n}\cdot\frac{1}{f\Big(\frac{u_0(X)}{\widetilde{D}_h x^n}\Big)}\cdot\frac{x^{n+1}-x^{n}}
%{\tau},x^{n+1}-x^{n}\right\rangle_{\mathcal{E}}   \leq 0. \nonumber
%&\geq \Big( d\big(\frac{f_{0}(X)}{D_{X}x^{n+1}}\big)+f_{0}(X)\frac{1-2x^{n}}{x^{n}(1-x^{n})}\Big|x^{n}-x^{n+1}\Big). \nonumber
%\end{align}
%\end{subequations}
%This completes the proof. \hfill$\Box$

%Finally  we provide the optimal rate convergence analysis for the schemes \eqref{equ:numnum}.
\begin{thm}%\label{PMELthm:convergence}
\label{convergence}
Assume that there exist positive constants $b_{u_0}$ and $B_{u_0}$, such that the initial condition $u_0(X)$ satisfies $b_{u_0} \leq u_0(X) \leq B_{u_0}$. Assume that $f(\cdot)$ and $H(\cdot)\in C^{2}([0, \infty))$, and $V_c(\cdot)$, $V_e(\cdot)$, $W_c(\cdot)$,  and $W_e(\cdot)\in C^2(\overline{\Omega})$. Denote by $x_e\in\Omega$ the exact solution to the trajectory equation \eqref{eqtra} with sufficient regularity and $x_{h}\in\mathcal{Q}$ the numerical solution to the numerical scheme~\eqref{equ:numnum}. The numerical error function is defined by%   based on the total energy $f\ln f$ or the numerical scheme based on the total energy $1/f$  \eqref{eqtranum1},
\begin{equation}\label{error-function-1}
 e_i^n = x_{e_i}^n - x_{i}^n ,
\end{equation}
where $x_{e}^n,\ x^n\in\mathcal{Q}$, $0 \le i \le M$, $n=0,\cdots,N$. %\mathbb{N}^+
Then
\begin{compactenum}
\item[$\bullet$]{
$e^{n}=(e^{n}_0,\cdots,e^{n}_M)$ satisfies}
$$\quad  \| e^{n}\|:=\langle e^{n},e^{n}\rangle_{\mathcal{E}} \le C (\tau + h^2).$$
%where    $T$ is the terminal time.
% \|_{2} \le e^{\frac12 C_0 T} \left( \tau  \sum_{k=1}^{n+1} \| \zeta^k \|_2^2 \right)^{1/2}
%$| \zeta^{n+1}_i | \le C (\tau + h^2)$, $1 \le i \le N-1$
\item[$\bullet$]{$\widetilde{D}_h e^{n}=(\widetilde{D}_h e^{n}_0,\cdots,\widetilde{D}_h e^{n}_M)$ satisfies}
$$\quad  \| \widetilde{D}_h e^{n} \| \leq C(\tau+h^2).$$
Moreover, numerical error between the numerical solution $u_h^{n}$ and exact solution $u_e^{n}$ to the  problem \eqref{equ:GFP} is estimated by %of equation \eqref{equ:origin}
$$\quad  \|u_h^{n}-u_e^{n}\| \leq C(\tau+h^2), ~n=0,\cdots,N,$$
where $C$ is a positive constant, $h$ is the grid spacing, and $\tau$ is the time step size.
\end{compactenum}
\end{thm}
The theorem can be proved by using a careful high-order asymptotic expansion and two-step error estimates, as shown in the work \cite{C.H.Duan(2019)}. The details of the proof are quite involved and will be presented in our future work.

The following theorem presents that our numerical method also preserves steady states.
\begin{thm}\label{steady} (Steady-state preserving)
Assume $u_0(x)>0$ and the discrete energy $E_N^n$ is bounded below. The fully discrete scheme \reff{equ:numnum} and \reff{equ:numbou} is steady-state preserving in the sense that, for fixed $h$,  as the time step $n\rightarrow \infty$, the numerical solution $\{x_i^{\infty}\}_{i=0}^{M}\in \mathcal{Q}$ is a numerical solution to a boundary-value problem (BVP) for steady states:
\begin{equation}\label{equ:steady}
\left\{
\begin{aligned}
&-\partial_X \left[\frac{u_0(X)}{\partial_X x^{\infty}}H'\Big(\frac{u_0(X)}{\partial_X x^{\infty}}\Big)-H\Big(\frac{u_0(X)}{\partial_X x^{\infty}}\Big)\right]
-u_0(X)V'(x^{\infty}) \\
&\qquad \qquad -u_0(X)\int_{\Omega}W'(x^{\infty}(X,t)-y^{\infty}(Y,t))u_0(Y)dY=0, \\
&x^{\infty} |_{\partial\Omega}= X|_{\partial\Omega}.
\end{aligned}
\right.
\end{equation}
\end{thm}
\noindent\textbf{Proof:}
%\begin{equation}\label{equ:steady}
%\begin{split}
%&\partial_X \left[\frac{u_0(X)}{\partial_X x^{\infty}}H'\Big(\frac{u_0(X)}{\partial_X x^{\infty}}\Big)-H\Big(\frac{u_0(X)}{\partial_X x^{\infty}}\Big)\right]\\
%&-u_0(X)V'(x^{\infty})-u_0(X)\int_{\Omega}W'(x^{\infty}(X,t)-y^{\infty}(Y,t))u_0(Y)dY=0.
%\end{split}
%\end{equation}
The energy dissipation in Lagrangian coordinate reads
\begin{equation}\label{equ:energylawL}
      \frac{d}{dt}E^{total}=-\Delta,
  \end{equation}
  where
%   the total energy $$E^{total}=\int_{\Omega}\left[H\Big(\frac{u_0(X)}{\partial_X x}\Big)\partial_X x+u_0(X)V(x(X,t))\right]dX+\frac{1}{2}\int_{\Omega}\int_{\Omega}W\left(x(X,t)-y(Y,t)\right)u_0(X)u_0(Y)dYdX,$$ and
  the entropy production
  \begin{equation*}
  \begin{split}
  \Delta&=\int_{\Omega} f\Big(\frac{u_0(X)}{\partial_X x}\Big)\cdot\frac{\partial_X x}{u^2_0(X)}\left|\partial_X \left[\frac{u_0(X)}{\partial_X x}H'\Big(\frac{u_0(X)}{\partial_X x}\Big)-H\Big(\frac{u_0(X)}{\partial_X x}\Big)\right]\right.\\
  &\qquad \qquad \left.+u_0(X)V'(x)+u_0(X)\int_{\Omega}W'(x(X,t)-y(Y,t))u_0(Y)dY\right|^2dX.
  \end{split}
  \end{equation*}
The energy dissipation \eqref{equ:energylawL} implies that steady states are achieved when the trajectory $x^{\infty}$ satisfies the BVP \reff{equ:steady}.
%\begin{equation}\label{equ:steady} \begin{split}
%&\partial_X \left[\frac{u_0(X)}{\partial_X x^{\infty}}H'\Big(\frac{u_0(X)}{\partial_X x^{\infty}}\Big)-H\Big(\frac{u_0(X)}{\partial_X x^{\infty}}\Big)\right]\\
%&-u_0(X)V'(x^{\infty})-u_0(X)\int_{\Omega}W'(x^{\infty}(X,t)-y^{\infty}(Y,t))u_0(Y)dY=0. \end{split}\end{equation}

Since the discrete total energy $E_N^n$ decreases monotonically and is bounded below, the limit  $E_N^{\infty}:=\lim\limits_{n\rightarrow \infty} E_N^n$ exists.  By the discrete energy dissipation law \eqref{eqi:ene_total}, we have $$\lim\limits_{n\rightarrow \infty} \left\langle\frac{u_0^2(X)}{\widetilde{D}_h x^n}\cdot\frac{1}{f\Big(\frac{u_0(X)}{\widetilde{D}_h x^n}\Big)}\cdot\frac{x^{n+1}-x^{n}}
{\tau},\frac{x^{n+1}-x^{n}}{\tau}\right\rangle_{\mathcal{E}}=0,\  \ x^{n+1},x^{n}\in\mathcal{Q}. $$
This along with the fully discrete scheme~\eqref{equ:numnum} deduces that the numerical solution $\{x_i^{\infty}\}_{i=0}^{M}$ solves the difference equations
\begin{equation}\label{equ:steadyNum}
\left\{
\begin{aligned}
&-d_h \left[\frac{u_0}{D_h x^{\infty}}\cdot H'\Big(\frac{u_0}{D_h x^{\infty} }\Big)-H\Big(\frac{u_0}{D_h x^{\infty}}\Big)\right]_i
-u_{0_i} V(x^{\infty}_{i})\\
 &\qquad \qquad -u_{0_i}\left\langle W'(x^{\infty}_{i}-y^{\infty}),u_0(Y)\right\rangle_{\mathcal{E}}=0,\ i=1,\cdots,M-1, \\
&x^{\infty}_0=X_0 ~ \mbox{and~} x^{\infty}_M=X_M,
\end{aligned}
\right.
\end{equation}
which are exactly a second-order finite difference discretization to the BVP \reff{equ:steady}. Thus, $\{x_i^{\infty}\}_{i=0}^{M}$ is a numerical solution to the BVP \reff{equ:steady}. $\hfill\Box$

\section{Numerical Results}
\label{sec:5}
In this section, we present numerical examples to demonstrate advantages of our numerical methods in solving the Fokker-Planck equation with nonlinear diffusion,  various nonlocal interaction kernels, and free boundaries with waiting time phenomena. In the following numerical simulations, we first solve the trajectory equation~\eqref{eqtra} with the initial and boundary conditions \eqref{eqtraFbou}-\eqref{eqtraFini} using the fully discrete scheme~\eqref{equ:numnum} with \eqref{eqnumbou},  and then approximate the density function $u$ in~\eqref{equ:conservationL} by \eqref{numdist}-\eqref{numdist2}.
%with one-well potential and double-well potential, waiting time phenomenon, nonlinear Fokker--Planck equations which describe the Boson gas, and aggregation-diffusion equations with the  smooth kernel and compacted support kernel.

We define the error of numerical solutions $e_h$ in $\mathcal{L}^2$, $\mathcal{L}^1$, and $\mathcal{L}^{\infty}$ norms as follows
\begin{equation}\label{L2_u}
\|e_h\|_2^2:=\left(\frac{1}{2}e_{h_0}^2 h_{x_{0}}+\sum\limits_{i=1}^{M-1}e_{h_i}^2 h_{x_i}+\frac{1}{2}e_{h_M}^2 h_{x_M}\right),
\end{equation}
\begin{equation}\label{L1}
\|e_h\|_1=\left(\frac{1}{2}|e_{h_0} |h_{x_{0}}+\sum\limits_{i=1}^{M-1}|e_{h_i}| h_{x_i}+\frac{1}{2}|e_{h_M}| h_{x_M}\right),
\end{equation}
and
\begin{equation}\label{Linf}
\|e_h\|_{\infty}=\max\limits_{0\leq i\leq M}\{|e_{h_i}|\},
\end{equation}
where $e_h=(e_{h_0},e_{h_1},\cdots,e_{h_M})$. Note that the error of the numerical trajectory $e^x_h:=x_e-x$ is given on the mesh with
\[
h_{x_i}=h,  \ \ 0\leq i\leq M,
\]
where $x_e$ is the exact solution of the trajectory on the grid and $h$ is a uniform grid spacing; and the error of the density $e^u_h:=u_e-u$ is given on the mesh with
\begin{equation*}
 h_{x_i}=\frac{x_{i+1}-x_{i-1}}{2}, \ \ 1\leq i \leq M-1, \ \ \
 h_{x_0}=x_{1}-x_{0}, \ \
 h_{x_M}=x_{M}-x_{M-1},
 \end{equation*}
where $u_e$ is the exact solution of density on the grid.
\subsection{Porous Medium Type of Equations}
In this section, we focus on the porous medium type of equations with one-well and double-well potentials, which have waiting time phenomena.
%\subsubsection{l}

{\bf Example 1: Porous medium equation with one-well potential }

We consider the following nonlinear equation
\begin{equation}\label{eqPME}
\partial_t u=\partial_x\left[u\partial_x\left(\frac{m}{m-1}u^{m-1}+\frac{x^2}{2}\right)\right], \ x\in[-2, 2],
\end{equation}
which corresponds to the FP equations with $f(u)=u$, $H(u)=\frac{1}{m-1}u^m$, $V(x)=\frac{x^2}{2}$, and $W=0$. This equation has been studied in \cite{J.Carrillo(2000),Z.Sun(2018)}. This equation with zero-flux boundary conditions has a steady state
%\begin{equation}\label{eqPMEini}
$$u_{\infty}(x)=\left(C-\frac{m-1}{2m}|x|^2\right)^{\frac{1}{m-1}},$$
where the constant $C$ is determined by ensuring the mass conservation.  Moreover, the relative total energy $E(t|\infty)=E(u(t))-E(u_{\infty})$ decays exponentially, i.e.,
$E(t|\infty)\leq E(0|\infty)e^{-2t}$. Here the decay rate $e^{-2t}$ is sharp.

Let $m=2$. We take the same compact supported initial data as in \cite{Z.Sun(2018)}:
\begin{equation}\label{InitCond}
u(x,0)=\max\{1-|x|,0\},\ x\in[-2, 2].
\end{equation}
The corresponding stationary solution is given in \cite{J.Carrillo(2000),Z.Sun(2018)}:
\begin{equation}\label{StationU}
u_{\infty}=\max\left\{\Big(\frac{3}{8}\Big)^{\frac{2}{3}}-\frac{x^2}{4},0\right\}.
\end{equation}

\begin{figure}[H]
\centering
\subfigure[Steady State]
{\includegraphics[scale=1.3]{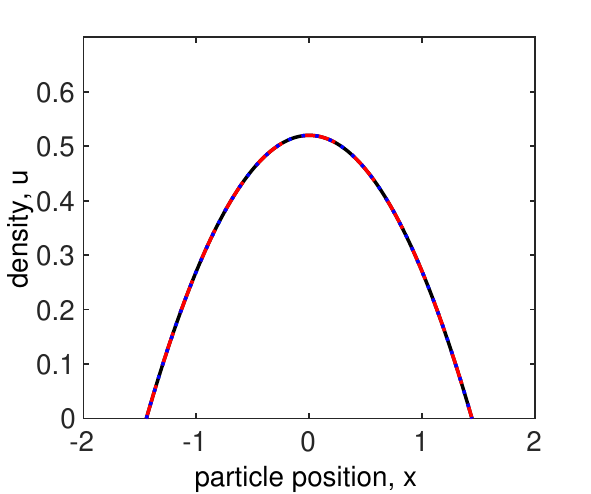}}
\subfigure[Zoomed-in plot near $x=3^{\frac{1}{2}}$]
{\includegraphics[scale=1.3]{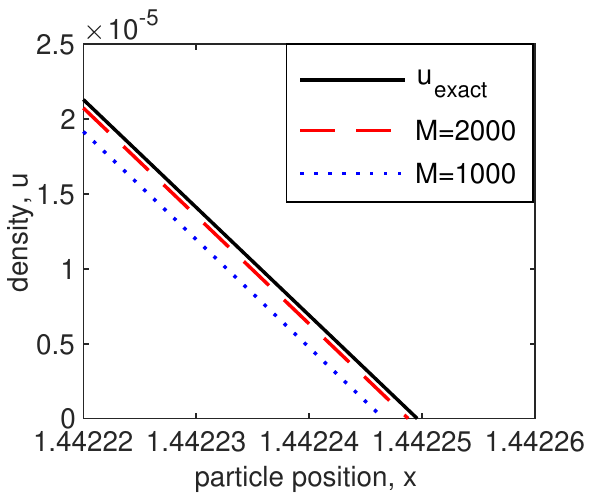}}
\caption{Numerical solution and the exact solution of the steady-state density of \textbf{Example 1} with $\tau=1/2000$ and $m=2$. }
\label{fig:PMEOneDen}
\end{figure}

We numerically simulate the problem on meshes with $M=1000$ and $M=2000$ up to time $T=10$, at which the problem can be identified as the steady state.  Fig. \ref{fig:PMEOneDen} (a) shows the numerical solution and exact solution \reff{StationU} at the steady state. Clearly, one can observe that the numerical solution can approximate the exact solution accurately. To have detailed comparison, we display a zoomed-in plot near the free boundary $x=3^{\frac{1}{2}}$ in Fig. \ref{fig:PMEOneDen} (b). As the mesh refines, the numerical solution can effectively approach the exact solution near the free boundary \emph{without any oscillation}.

%\begin{figure}[H]
%\centering
%\subfigure[Relative energy $E(t|\infty)$]
%{\includegraphics{fig/OneWellEne.eps}}%[scale=1.2]
%\subfigure[\scriptsize $\frac{E(t|\infty)}{E(0|\infty)}$]
%{\includegraphics{fig/OneWellReEne.eps}}%ReEneShu.png [scale=1.2]
%\caption{The evolution of the relative energy of \textbf{Example 1} with $h=\tau$ and $m=2$.}
%\label{fig:PMEOneEne}
%\end{figure}

\begin{figure}[H]
\centering
\subfigure[Relative energy $E(t|\infty)$]
{\includegraphics[scale=1.1]{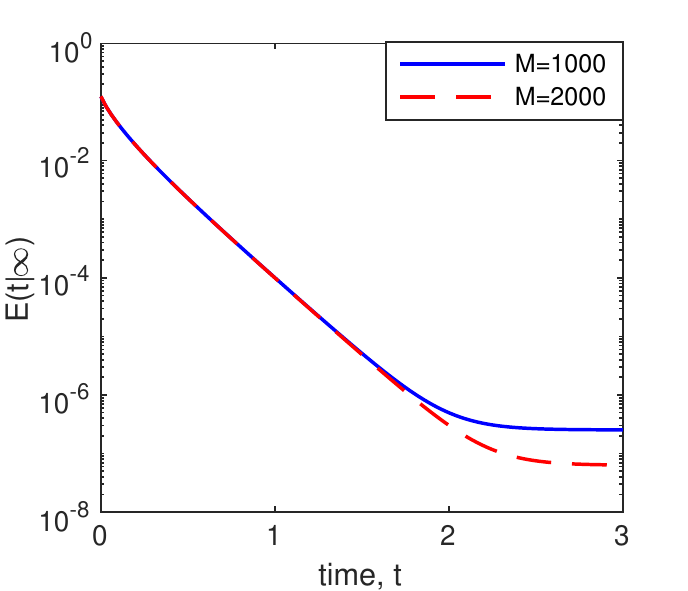}}
\subfigure[\scriptsize $\frac{E(t|\infty)}{E(0|\infty)}$]
{\includegraphics[scale=1.1]{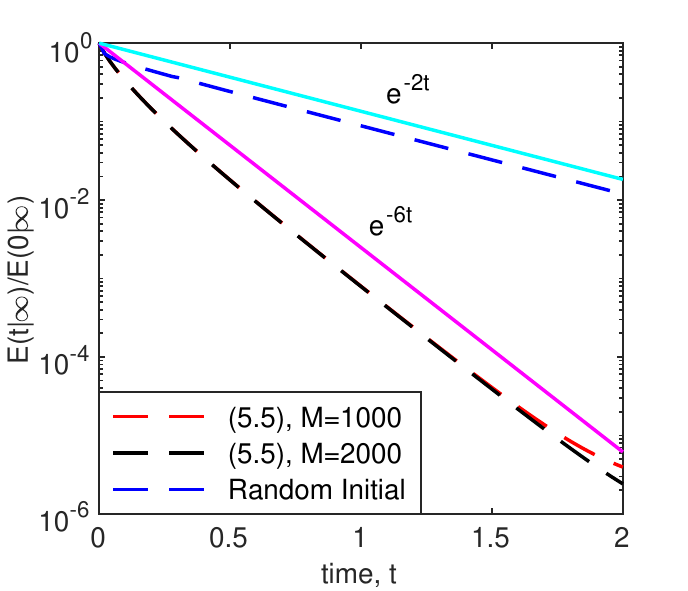}}
\caption{The evolution of the relative energy of the equation \reff{eqPME} with the initial condition \reff{InitCond} and a random initial condition in \textbf{Example 1}  ($\tau=h/2$ and $m=2$).}
\label{fig:PMEOneEneln}
\end{figure}

%\begin{figure}[H]
%\centering
%\subfigure[The random initial function $u_0$]
%{\includegraphics[scale=1.1]{fig/OneWellU0.eps}}
%\subfigure[\scriptsize $\frac{E(t|\infty)}{E(0|\infty)}$ ]
%{\includegraphics[scale=1.1]{fig/OneWellRandomEne.eps}}
%\caption{The evolution of the relative energy  for the random initial function of \textbf{Example 1}  ($M=1000$, $\tau=1/1000$ and $m=2$).}
%\label{fig:PMEOneEneRandom}
%\end{figure}

Furthermore, we show in Fig. \ref{fig:PMEOneEneln} (a) the relative total energy. One can see that the relative total energy decays monotonically and remains positive as time evolves for both $M=1000$ and $M=2000$. To further check the decay rate, we present the relative total energy rescaled by its initial value in Fig. \ref{fig:PMEOneEneln} (b), which displays that the decay rate is around  $e^{-6t}$ for the equation \reff{eqPME} with the initial condition \reff{InitCond}.  To investigate more on the decay rate, we also study the equation \reff{eqPME} with a random initial condition.  The dashed blue curve shown in Fig. \ref{fig:PMEOneEneln} (b) demonstrates that $\frac{E(t|\infty)}{ E(0|\infty)} \leq e^{-2t}$,  being consistent with the theoretical result \cite{J.Carrillo(2000)}.

%\begin{figure}
%%\begin{minipage}[t]{0.5\linewidth}
%\captionsetup{font={scriptsize}}
%\centering
%\subfigure[\scriptsize Steady State at time $t=10$]
%{\includegraphics[width=5cm,height=4cm]{fig/OneWellDensity.eps}}\hspace{.3in}
%\subfigure[\scriptsize Zoomed in figure near $x=3^{\frac{1}{2}}$] %=\frac{1}{5}(2+6x+\frac{\pi}{2}\sin(2\pi x))
%{\includegraphics[width=5cm,height=4cm]{fig/OneWellDensityDe.eps}}
%\caption{ Numerical density and the exact solution at steady state of \textbf{Example 1}  ($m=2$)}
%\label{fig:PMEOneDen}
%\end{figure}

%\begin{figure}
%%\begin{minipage}[t]{0.5\linewidth}
%\captionsetup{font={scriptsize}}
%\centering
%\subfigure[\scriptsize  Relative energy $E(t|\infty)$]
%{\includegraphics[width=5cm,height=4cm]{fig/OneWellEne.eps}}\hspace{.3in}
%\subfigure[\scriptsize $\frac{E(t|\infty)}{E(0|\infty)}$] %=\frac{1}{5}(2+6x+\frac{\pi}{2}\sin(2\pi x))
%{\includegraphics[width=5cm,height=4cm]{fig/OneWellReEne.eps}}%ReEneShu.png
%\caption{The evolution of relative energy of  \textbf{Example 1} ($h=1/1000$, $\tau=1/1000$,$m=2$)}
%\label{fig:PMEOneEne}
%\end{figure}

\begin{table}[H]
\centering
\begin{tabular}{cccccccc}
\hline

\hline
  $h$    &$\tau$ &$\|e^u_h\|_2$    & Order  &$\|e^u_h\|_1 $    & Order & $\|e^u_h\|_{\infty}$    & Order \\ \hline
1/50&1/100 &1.015e-03&&2.894e-04& &5.193e-03&  \\ \hline
1/100 &1/400 &3.597e-04&1.497&7.292e-05&1.989&2.598e-03&0.999 \\  \hline
1/200 &1/1600 &1.273e-04&1.498&1.830e-05&1.994&1.300e-03&0.999 \\ \hline
1/400 &1/6400 &4.505e-05&1.499&4.584e-06&1.997&6.499e-04&1.000 \\
 \hline

 \hline
\end{tabular}
\caption{Numerical error and convergence order of the numerical solution  at time $T=10$ for {\bf Example 1}.}\label{table:Onefinal}
\end{table}

%\begin{table}[H]
%\centering
%\begin{tabular}{cccccccc}
%\hline
%
%\hline
% $h$    &$\tau$ &$ \mathcal{L}^2$-error $ (u) $    & Order  &$ \mathcal{L}^1$-error $ (u) $    & Order &$ \mathcal{L}^{\infty}$-error$ (u) $    & Order \\ \hline
%0.1&0.1 &3.987e-03&&1.768e-03& &1.2952e-02& \\ \hline
%0.2 &0.2 &1.417e-03&1.492&4.505e-04&1.973&6.488e-03&0.997\\ \hline
%0.4 &0.4 &5.025e-04&1.496&1.137e-04&1.986&3.247e-03&0.999\\ \hline
%0.8 &0.8 &1.779e-04&1.498&2.857e-05&1.993&1.624e-03&0.999\\
%
% \hline
%
% \hline
%\end{tabular}
%\caption{Numerical error and convergence rate of the numerical solution $u$ at time $T=10$ for {\bf Example 1}.}\label{table:Onefinal}
%\end{table}
We further consider the numerical accuracy of our numerical method with various mesh resolution. Table \ref{table:Onefinal} shows that the numerical error and convergence rate for the solution $u$ at time $T=10$ in the $\mathcal{L}^1$, $\mathcal{L}^2$  and  $\mathcal{L}^{\infty}$ norms. We observe that our numerical method has convergence order of $2$ in the $\mathcal{L}^1$ norm, $\frac{3}{2}$ in the $\mathcal{L}^2$ norm, and $1$ in the $\mathcal{L}^{\infty}$ norm. It is believed that the low regularity near the free boundary accounts for the decrease of convergence order. We remark that, to the best of our knowledge, the convergence order of numerical schemes for the Fokker-Planck equations with the presence of free boundaries has not been studied in the literature.  %The convergence rate has not yet exhibited  at other recent works.

{\bf Example 2: Porous medium equation with a double-well potential}

In this example, we consider the degenerate Fokker--Planck equations with a double-well potential:
\begin{equation}\label{eqPMEDW}
\partial_t u=\partial_x\left\{u\partial_x\left[\nu u^{m-1}+V(x)\right]\right\}, \ x\in[-2, 2],
\end{equation}
in which $f(u)=u$, $H'(u)=\nu u^{m-1}$, $V(x)=x^4/4-x^2/2$, and  $W=0$.
The steady state is given by
\begin{equation}\label{equ:DWe_steady}
u_{\infty}=\left(\frac{c(x)-V(x)}{\nu}\right)_+^{1/(m-1)},
\end{equation}
where $c(x)$ is a piecewise constant \cite{J.Carrillo(2015), H.Liu(2016),Y.Qian(2019)}.

To demonstrate the accuracy of our
numerical schemes, we solve the problem \eqref{eqPMEDW} with a positive initial condition
\begin{equation}\label{equ:DWe}
u_0=\frac{M}{\sqrt{2\pi\sigma^2}}e^{-\frac{x^2}{2\sigma^2}},\ x\in[-2,2],
\end{equation}
where $M=4.2517\times 10^{-2}$ and $\sigma=\sqrt{0.2}$. Then $c(x)=-\frac{3}{16}$  in \eqref{equ:DWe_steady}.

%, at which the problem can be identified as the steady state.

\begin{figure}[H]
\centering
\includegraphics[scale=0.9]{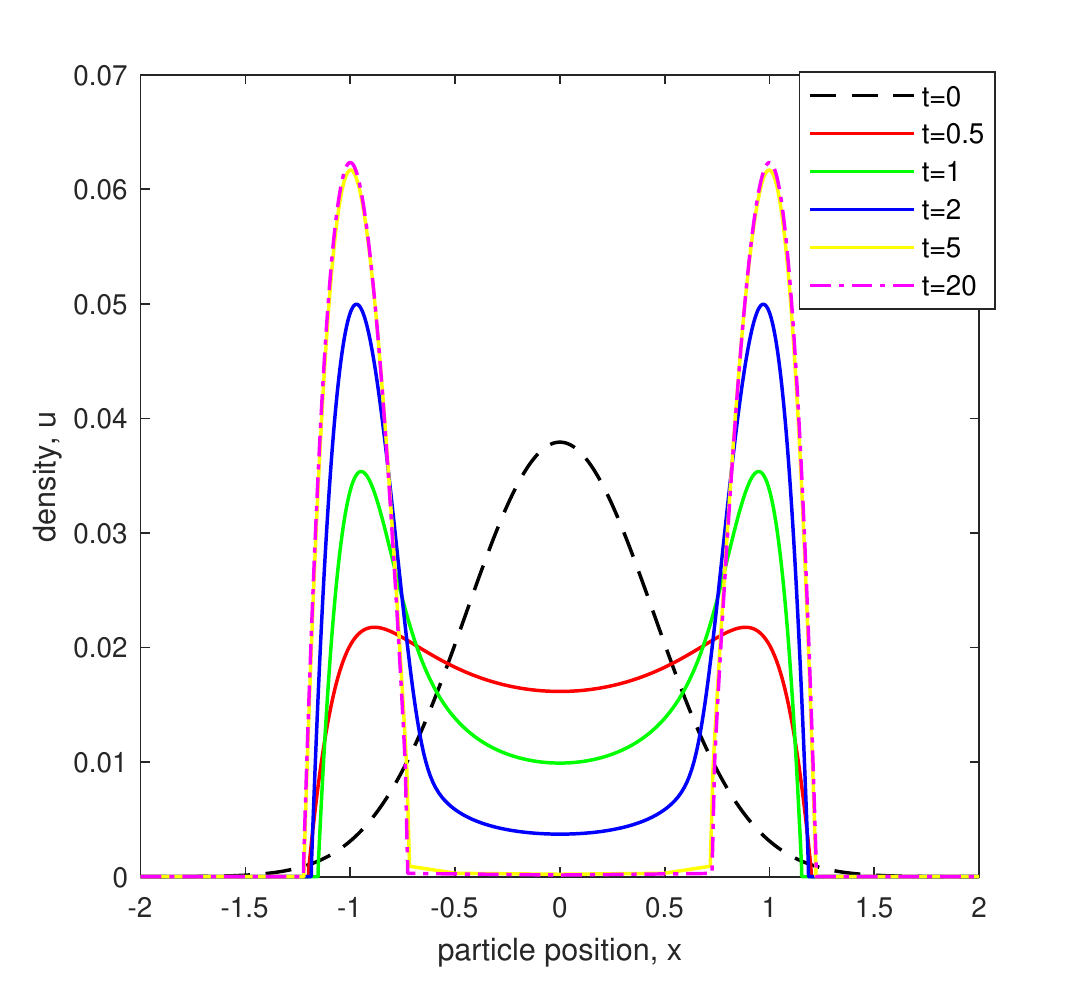}
\caption{The evolution of the numerical density in {\bf Example 2} with $h=1/1000$, $\tau=1/1000$, and $m=2$.}
\label{fig:DWEnergy}
\end{figure}

\begin{figure}[H]
\centering
\subfigure[Total Energy]
{\includegraphics[scale=1]{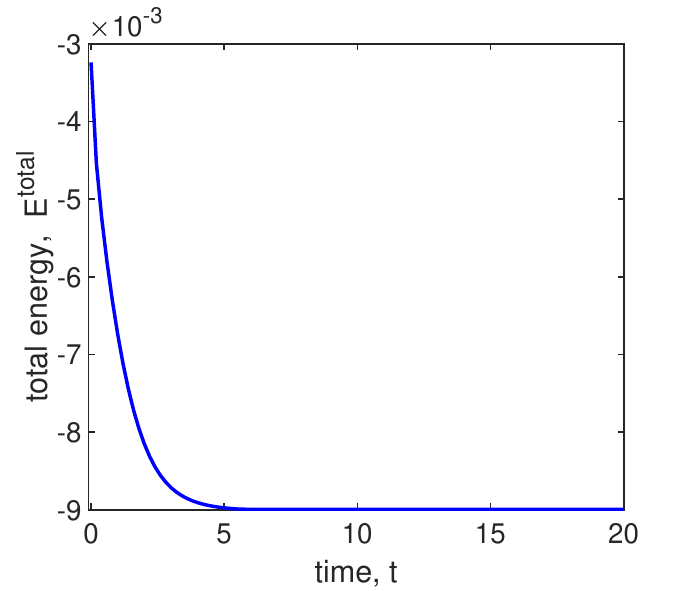}}
\subfigure[Particle positions]
{\includegraphics[scale=1.2]{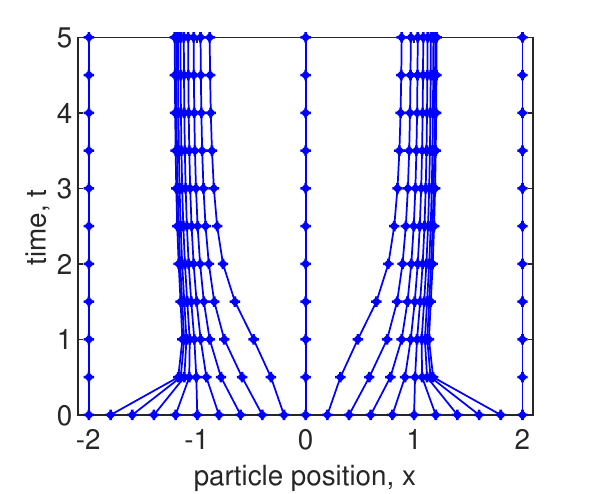}}
\caption{The evolution of the total energy and particle positions in {\bf Example 2} with $h=1/1000$, $\tau=1/1000$, and $m=2$.}
\label{fig:PoEnergy}
\end{figure}

\begin{figure}[H]
\centering
\includegraphics[scale=1.]{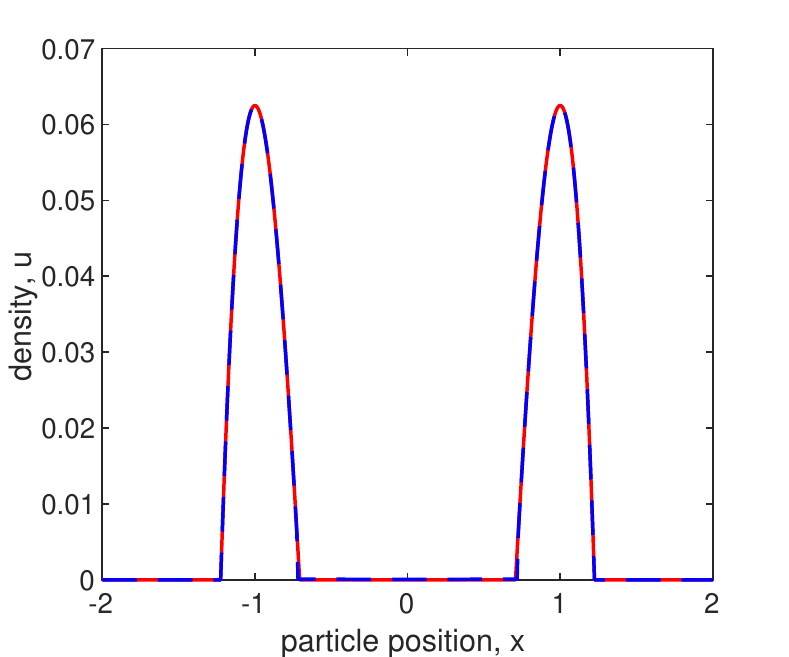}
\caption{The numerical solution $u_h$ (the blue line) and exact solution $u_e$ (the red dotted line) at the steady state  ($T=20$) in {\bf Example 2} with $h=1/1000$, $\tau=1/1000$, and $m=2$.}
\label{fig:DoubleSteady}
\end{figure}

Fig. \ref{fig:DWEnergy} shows the evolution of density $u$ with $h=1/1000$, $\tau=1/1000$, and $m=2$ up to time $T=20$. Due to the external double-well potential, the solution gradually splits into two parts, localizing at the two centers of the wells. Long-time simulation up to the time $T=20$ shows the asymptotic convergence of the numerical solution towards a steady state.  Fig. \ref{fig:PoEnergy} (a) and (b) show  the decay of the total energy and motion of particles that move towards the two centers with a finite speed. Fig. \ref{fig:DoubleSteady}  shows that the numerical solution and   exact solution at the steady state are almost identical.

\begin{table}[H]
\centering
\begin{tabular}{cccccc}
\hline

\hline
$h$    &$\tau$ &$ \|e^u_h\|_2$    & Order   &$\|e^u_h\|_{\infty}$   & Order  \\\hline
0.02&0.02 &5.606e-05& &5.962e-05& \\\hline
0.01 &0.005 &1.422e-05&1.979 &1.509e-05&1.983\\\hline
0.005 &0.00125 &4.170e-06&1.769 &3.781e-06&1.996\\
 \hline
$h$    &$\tau$ &$ \|e^x_h\|_2$    & Order   &$\|e^x_h\|_{\infty}$    & Order \\\hline
0.02&0.02 &1.619e-02 &&2.390e-02&\\\hline
0.01 &0.005 &4.334e-03&1.979&6.485e-03&1.882\\\hline
0.005 &0.00125 &1.106e-03&1.770&1.694e-03&1.936\\
 \hline

  \hline
\end{tabular}
\caption{Numerical error and convergence order of the numerical trajectory $x$ and density $u$ at time $T=0.1$ in {\bf Example 2}.}\label{table:DW1}
\end{table}
\begin{table}[H]
\centering
\begin{tabular}{cccccc}
\hline

\hline
 $h$    &$\tau$ &$\|e^u_h\|_2$    & Order   &$\|e^u_h\|_{\infty}$    & Order  \\\hline
 1/50&1/50 &1.3554e-02&&1.5146e-02& \\\hline
 1/100 &1/200 &9.7065e-03&0.4817&1.1056e-02&0.4541\\\hline
 1/200 &1/800 &6.9297e-03&0.4862&7.9973e-03&0.4673\\\hline
 1/400 &1/3200 &4.9353e-03&0.4897&5.7471e-03&0.4767\\
 \hline

 \hline
\end{tabular}
\caption{Numerical error and convergence order of the numerical density $u$ at the steady state in {\bf Example 2}.}\label{table:DW2}
\end{table}

Table \ref{table:DW1} presents the numerical error and convergence order of both the numerical trajectory $x$ and density $u$ at time $T=0.1$. Note that the reference exact solution is obtained numerically on a rather refined mesh with $h=10^{-5}$ and $\tau =10^{-6}$. One observes that the numerical method is roughly second order accurate in space and first order accurate in time at time $T=0.1$. However, the numerical convergence order degenerates to $0.5$ in space at the steady state in Table \ref{table:DW2}, due to the lower regularity of the solution close to  the free boundaries; cf. the solution profile with a support  in Fig. \ref{fig:DWEnergy}.

{\bf Example 3: Waiting time phenomena with nonlocal interactions}

In this example, we consider a nonlinear diffusion system with  the waiting time phenomena, which are commonly found in the porous medium type of diffusion \cite{D.G.Aronson(1983),J.L.Vazquez(2007),C.H.Duan(2018)}.
The challenge lies in establishing an efficient algorithm for the calculation of waiting time. One related algorithm has been proposed in the work \cite{C.H.Duan(2018)} to compute the waiting time for the porous medium equation. We now extend the algorithm to solve a more complex problem~\eqref{equ:GFP} with $f(u)=u$, $H(u)=\frac{\nu}{m}u^m$, $V(x)=-\frac{(x+\frac{\pi}{2})^2}{2}$, $W(x)=\frac{4-\theta}{8}|x|$, and the initial condition
$$u_0(x)=\left\{\frac{m-1}{m}[(1-\theta)\sin^2(x)+\theta \sin^4(x)]\right\}^{1/(m-1)},\ x\in[-\pi,0],$$
where  $\theta\in[0,1]$ is a parameter. Note that the nonlocal interaction kernel $W(x)$ here is not differentiable at $x=0$. To address this issue, we split the integral \reff{eqtraSS} into two integrals on two domains with the non-differentiable location as the integration upper and lower limits of such two integrals. In each integral, the $W'(0)$ is understood as $W'(0\pm)$.  We employ a mesh with a total spatial grid number $M=100$ and a time step size $\tau=1/100$.
Let $m=2$ and $\nu=2$.
%Now we  success the algorithm in Sec. 3.3 and present the detail  algorithm as following:
%To compute the waiting time, we propose the following algorithm.
%
%{\bf Algorithm for Waiting time}
%\begin{itemize}
%\item {\bf Step 1.}\ \ For time $t^n, n=0,1, \cdots$, solve the trajectory equation \eqref{eqtra} with \eqref{eqtraFini} and \eqref{eqtraFbou} using the fully discrete scheme  \eqref{equ:numnum} with the boundary conditions $\partial_t x = 0$ replaced by $x^{n+1}_0 = \xi_1^0$ and $ x^{n+1}_M = \xi_2^0$, where $\xi_1^0$ and $\xi_2^0$ are defined in  \eqref{equ:boundaryleft} and \eqref{equ:boundaryright}, respectively.
%
%
%
%
%    Check the criterion \eqref{WTC2_1} for $x^{n+1}$. If it is not valid, goto next time step. If it is valid, then set $t^*_h = t^{n+1}$. $n^*= n+1$ and  goto Step 2.
%\item {\bf Step 2.}\ \ For time $t^n, M=n^*,n^*+1, \cdots$, solve the trajectory equation  \eqref{eqtra} with \eqref{eqtraFini} and \eqref{eqtraFbou} by the fully discrete scheme \eqref{equ:numnum} with \eqref{eqnumbou}.
%\end{itemize}
\begin{figure}
\centering
\subfigure[\scriptsize $t=0$]
{\includegraphics[scale=0.85]{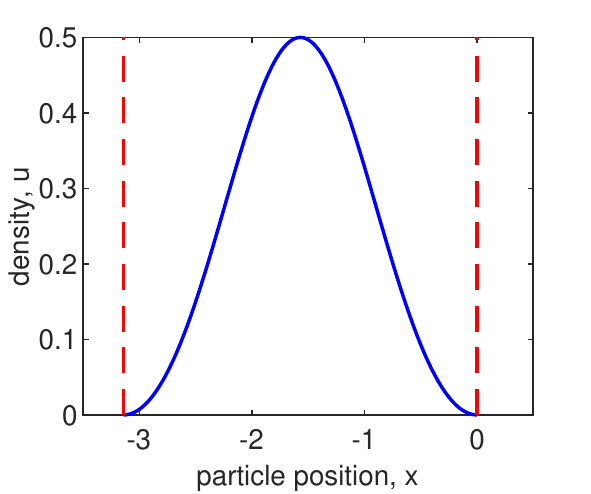}}\hspace{-.1in}
\subfigure[\scriptsize $t=0.1$]
{\includegraphics[scale=0.85]{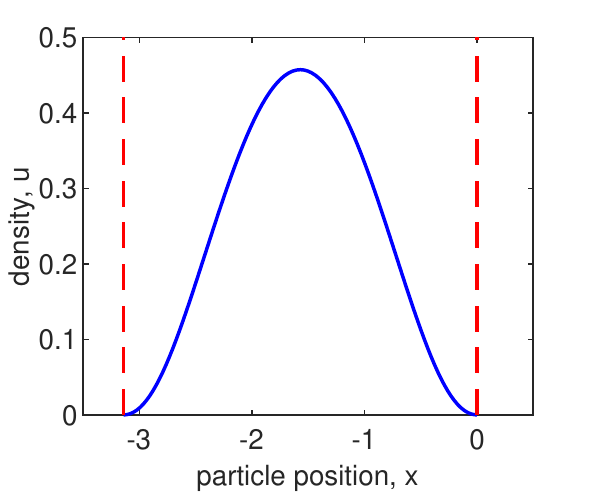}}\hspace{-.1in}
\subfigure[\scriptsize $t=0.31$]
{\includegraphics[scale=0.85]{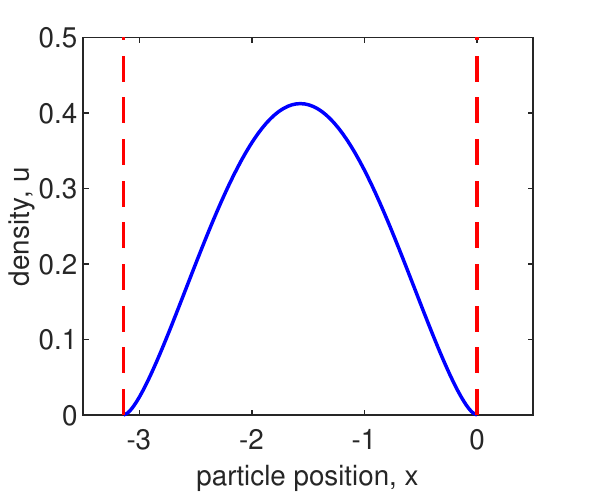}}\hspace{-1in}\\
\subfigure[\scriptsize $t=0.5$]
{\includegraphics[scale=0.85]{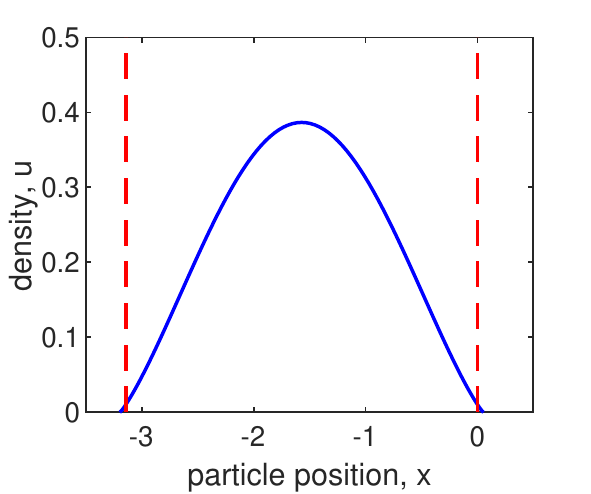}}\hspace{-.1in}
\subfigure[\scriptsize $t=0.8$]
{\includegraphics[scale=0.85]{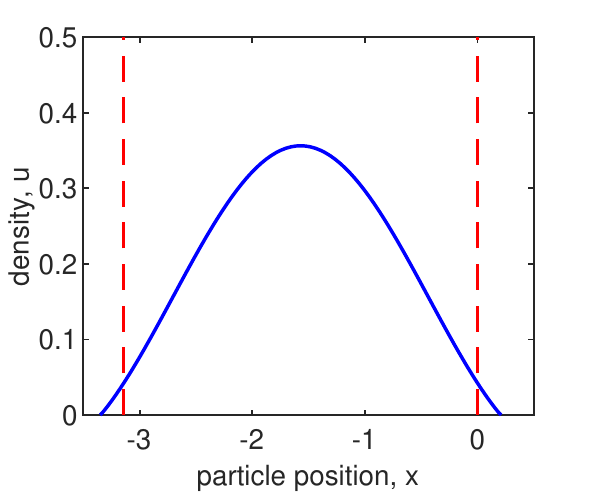}}\hspace{-.1in}
\subfigure[\scriptsize $t=1$]
{\includegraphics[scale=0.85]{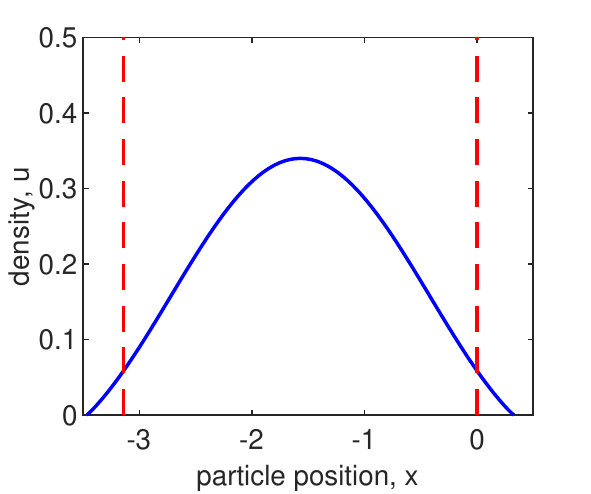}}\hspace{-.1in}%\\
\caption{Evolution of the profile of the density $u$ from $t=0$ to $t=1$ with $\theta=0.25$ in \textbf{Example 3}.}
\label{fig:WTDensity}
\end{figure}

\begin{figure}
\centering
\subfigure[]{\includegraphics[scale=1.3]{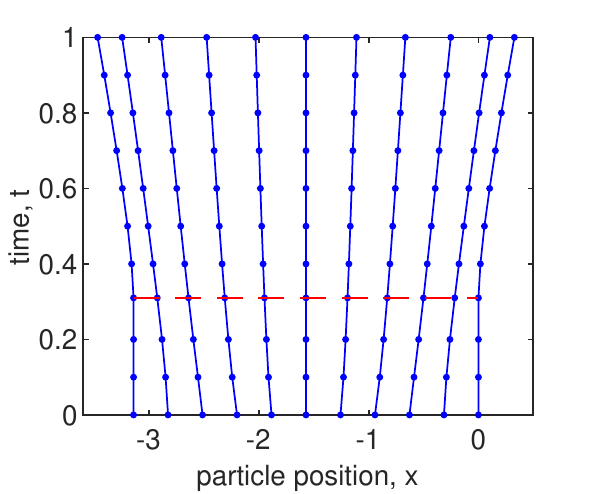}}\hspace{-.1in}
\subfigure[]{\includegraphics[scale=1.3]{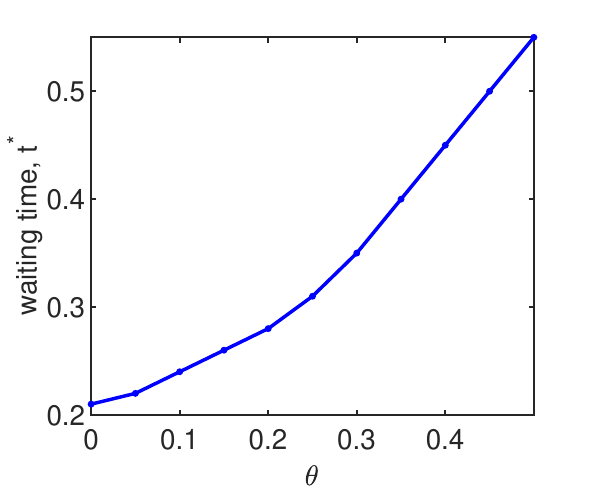}}
\caption{(a) Evolution of particles with $\theta=0.25$; (b) The relationship between $\theta$ and the waiting time $t^*$ in \textbf{Example 3}.}
\label{fig:WTtheta}
\end{figure}

%\begin{table}[H]
%\centering
%\begin{tabular}{cccccc}
%\hline
%
%\hline
%$M$    &$\tau$ &$\|e_h^u\|_2 $    & Order   &$\|e_h^u\|_{\infty} $  & Order  \\\hline
% 100&1/100 &9.637e-04&&1.046e-03& \\\hline
% 200 &1/400 &2.460e-04&1.970&2.664e-04&1.974\\ \hline
% 400 & 1/1600 &6.121e-05 &2.007&6.621e-05&2.008 \\ \hline
%$M$    &$\tau$& $\|e_h^x\|_2 $     & Order   &$\|e_h^x\|_{\infty}  $     & Order\\\hline
% 100&1/100 &1.020e-03&&2.316e-03&\\ \hline
% 200 &1/400 &2.991e-04&1.970&1.153e-03&1.006 \\ \hline
%400&1/1600 &9.165e-05&2.007&5.761e-04&1.001\\ \hline
%
%  \hline
%\end{tabular}
%\caption{Numerical error and convergence order of the numerical trajectory $x$ and density $u$ at time $t=0.1$ in {\bf Example 3}.}\label{table:waitingtime}
%\end{table}

%\begin{figure}
%\begin{minipage}[t]{0.5\linewidth}
%%\captionsetup{font={scriptsize}}
%\centering
%{\includegraphics[width=5cm,height=4cm]{fig/WTParticle_new.eps}}\hspace{.3in}
%\caption{Evolution of particles in \textbf{Example 3}.}%($h=1/100$, $\tau=1/10000$, $m=2$)
%%\subfigure[\scriptsize Zoomed in figure near $t=0$] %=\frac{1}{5}(2+6x+\frac{\pi}{2}\sin(2\pi x))
%\label{fig:PMEOneParWT}
%\end{minipage}
%
%\begin{minipage}[t]{0.5\linewidth}
%    \captionsetup{font={scriptsize}}
%     \centering
%\includegraphics[width=5cm,height=4cm]{fig/WT_theta.eps}
%\caption{ The relationship between $\theta$ and the waiting time $t^*$ in \textbf{Example 3}} % ($h=1/1000$, $\tau=1/2000$, $m=2$)  for \textbf{Example 3}
%\label{fig:PMEWT_theta}
% \end{minipage}
%\end{figure}

Fig. \ref{fig:WTDensity} shows evolution of the profile of the density $u$ with $m=2$ from $t=0$ to $t=1$. The estimated waiting time is about $t^*=0.31$.  During the time evolution, the numerical solution is free of oscillation near free boundaries and the boundaries have finite-speed propagation after exceeding the waiting time.  Fig. \ref{fig:WTtheta} (a) shows trajectories of particles as time evolves. One can find that  the  boundaries remain stationary up to time $t^*=0.31$, as indicated by the red dash line, and then move outward at a finite speed. Fig. \ref{fig:WTtheta} (b) displays the relationship between the waiting time $t^*$ and $\theta$. We can see that the waiting time $t^*$ increases monotonically as $\theta$ grows.

%Table \ref{table:waitingtime} shows the convergence order for $x$ and $u$ in $\mathcal{L}^2$ and $\mathcal{L}^{\infty}$ norms. One observes that the convergence order is about $2$ in space and $1$ in time, being consistent with the Theorem \ref{convergence}, although the nonlocal interaction kernel $W(x)$ does not meet the assumptions.

%The order decreases to 1st   in time and space, mainly because the regularity becomes low before the waiting time $t^*$.

%\begin{figure}
%%\captionsetup{font={scriptsize}}
%\centering
%\subfigure[\scriptsize $t=0$]
%{\includegraphics[width=4cm,height=3cm]{fig/WTDensity_t=0.eps}}\hspace{.3in}
%\subfigure[\scriptsize $t=0.1$]
%{\includegraphics[width=4cm,height=3cm]{fig/WTDensity_t=0_1.eps}}\hspace{.3in}
%\subfigure[\scriptsize $t=0.315$]
%{\includegraphics[width=4cm,height=3cm]{fig/WTDensity_t=0_315.eps}}\hspace{1in}\\
%\subfigure[\scriptsize $t=0.5$]
%{\includegraphics[width=4cm,height=3cm]{fig/WTDensity_t=0_5.eps}}\hspace{.3in}
%\subfigure[\scriptsize $t=1$]
%{\includegraphics[width=4cm,height=3cm]{fig/WTDensity_t=1.eps}}\hspace{.3in}
%\subfigure[\scriptsize $t=2$]
%{\includegraphics[width=4cm,height=3cm]{fig/WTDensity_t=2.eps}}\hspace{1in}%\\
%\caption{ Waiting time: Evolution of density $u$ ($m=2$, $M=100$,  $\tau=1/10000$)  for \textbf{Example 3}}
%\label{fig:WTDensity}
%\end{figure}

\subsection{Nonlinear Fokker--Planck Equations}

{\bf Example 4: Generalized Fokker--Planck Equations for Boson Gas}

In this example, we focus on the generalized Fokker--Planck equation with a superlinear drift:
\begin{equation}
\partial_t u=\partial_x\left[xu(1+u^K)+\partial_x u\right],
\end{equation}
where $K$ is a positive constant. For $K>2$, the system exhibits a critical mass phenomenon \cite{N.B.Abdallah(2011)}, i.e., an initial distribution with supercritical mass leads to singularity at the origin.  The phenomenon has been confirmed numerically in \cite{M.Bessemoulin-Chatard(2012), H.Liu(2016),Z.Sun(2018)}, in which the numerical solution approximates singularity with precision dependent on mesh resolution, i.e., $\mathcal{O}(1/h)$ with $h$ being the grid spacing. In this work, we shall show that the approximation precision of singularity can be enhanced significantly by our numerical methods.

We take $f(u)=u(1+u^3)$, $H'(u)=\log(\frac{u}{\sqrt[3]{1+u^3}})$, $V(x)=\frac{\beta}{2}x^2$, and $W=0$ in  \eqref{equ:GFP}. The initial data is given by
$$u_0(x)=\frac{M}{2\sqrt{2\pi}}\Big(e^{-\frac{(x-2)^2}{2}}+e^{-\frac{(x+2)^2}{2}}\Big),\ x\in[-6, 6].$$
Thus the trajectory equation becomes
\begin{equation}\label{equ:BG}
\partial_t x=-\left[1+\Big(\frac{u_0(X)}{\partial_X x}\Big)^{3}\right]\cdot\left[\frac{1}{\partial_X x}\log\left(\frac{u_0(X)/\partial_X x}{\sqrt[3]{1+(u_0(X)/\partial_X x)^3}}\right)+\beta x\right].
\end{equation}

%{\color{brown}{The motion of particles is primarily determined by \eqref{equ:BG}. In a supercritical case, the right term may tend to infinity as particles accumulate at certain position, i.e., the solution blows up at a finite time.}}

%at a point $x_{i_{cen}}$, i.e. $\exists$  $0\leq i_0\leq M$, such that $x_{i_{cen}}=x_{i_0}$.   \zhou{Do not understand here.}

\begin{rem}
%\emph{Criteria for particles meet.}
%In the discrete scheme  \eqref{equ:numnum} , we find the solution $x^{n+1}\in\mathcal{Q}$, i.e., $x_{i}<x_{i+1}$, $i=0,\cdots,M-1$, so that the particles never touch each other.  Hence \eqref{numdist} - \eqref{numdist2} make sense.
When the distances between particles get less than machine precision, they are indistinguishable in numerical calculations and the numerical solution may blow up. To avoid losing accuracy of density $u$ in \eqref{numdist}, we merge those particles together in our numerical treatment and regard them as one particle with the density changed accordingly. As a rule of thumb, we choose a criterion with a tolerance $\varepsilon_0=10^{-9}$ and define
\begin{equation}\label{eq:blowupB}
\mathcal{R}:=\{x_i^{n+1}|x_{i+1}^{n+1}-x_{i}^{n+1}\leq\varepsilon_0,\ i=0,\cdots, M-1\}.
\end{equation}
If $\mathcal{R}\neq\emptyset$,  there must be some particles that have been merged together at time $t^{n+1}$. See the work \cite{C.H.Duan(2017)} for more details on the numerical implementation.
\end{rem}

\begin{figure}
\centering
\subfigure[Entropy]
{\includegraphics[scale=1.3]{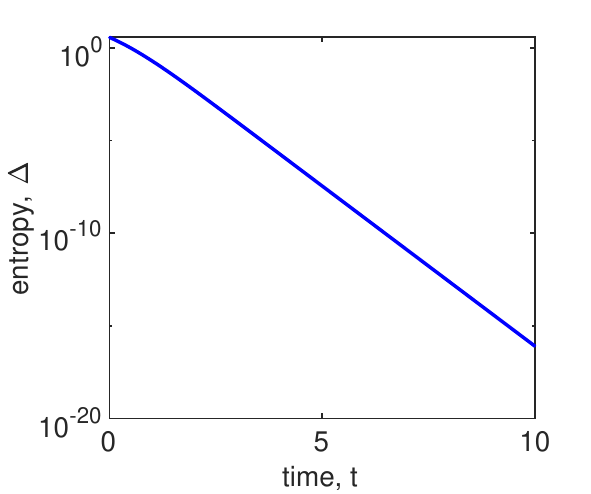}}
\subfigure[Particle trajectories] %=\frac{1}{5}(2+6x+\frac{\pi}{2}\sin(2\pi x))
{\includegraphics[scale=1.3]{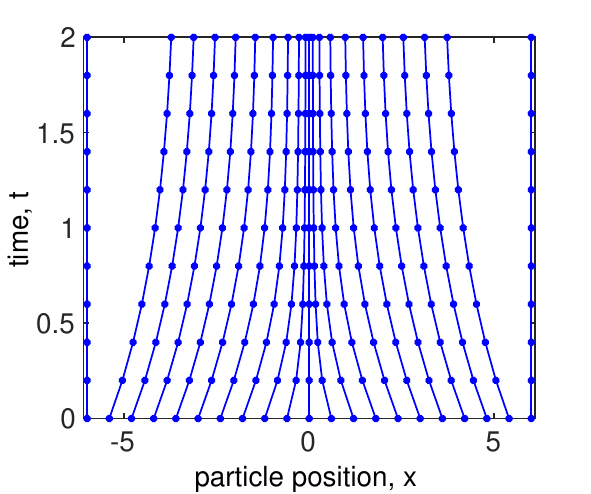}}
\caption{The evolution of entropy and particle trajectories for the case with subcritical mass $M=1$  and $\beta=1$ in  \textbf{Example 4} ($h=1/1000$, $\tau=1/1000$, $m=2$).}
\label{fig:BGEnergy}
\end{figure}

\begin{figure}
%\captionsetup{font={scriptsize}}
\centering
{\includegraphics[scale=1.1]{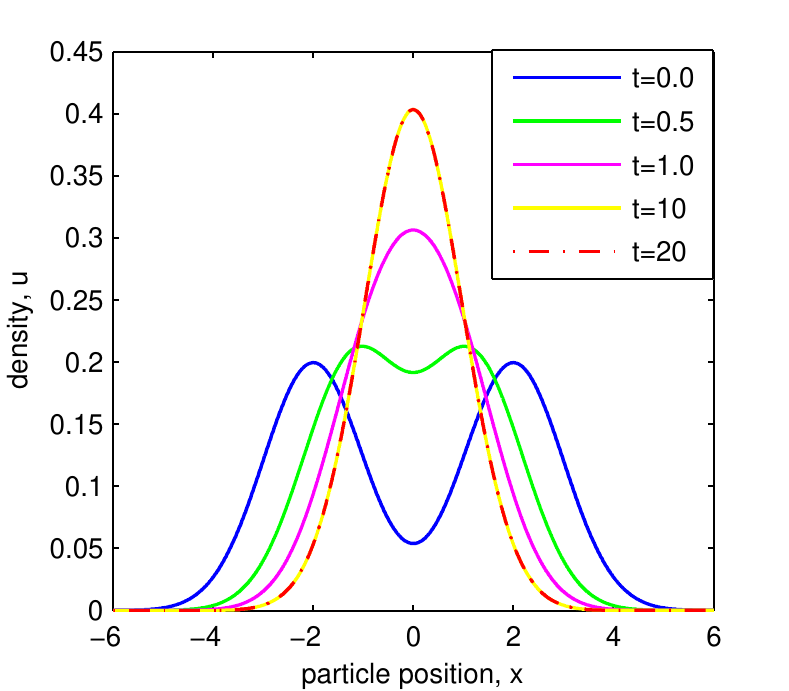}}
\caption{Evolution of  density $u$ with subcritical mass $M=1$ and $\beta=1$  ($h=1/1000$, $\tau=1/1000$, $m=2$) in \textbf{Example 4}.}
\label{fig:BGDensity}
\end{figure}

%\begin{figure}
%%\begin{minipage}[t]{0.5\linewidth}
%\captionsetup{font={scriptsize}}
%\centering
%\subfigure[\scriptsize Entropy]
%{\includegraphics[width=5cm,height=4cm]{fig/BGEntropyM=1.eps}}\hspace{.3in}
%\subfigure[\scriptsize Particle] %=\frac{1}{5}(2+6x+\frac{\pi}{2}\sin(2\pi x))
%{\includegraphics[width=5cm,height=4cm]{fig/BGParticleM=1.eps}}
%\caption{The evolution of entropy and particles  for subcritical mass $M=1$ for \textbf{Example 4} ($h=1/1000$, $\tau=1/1000$, $m=2$)}
%\label{fig:BGEnergy}
%\end{figure}

\begin{figure}
\centering
\subfigure[Time $t=0$ to $1$]
{\includegraphics[scale=1.1]{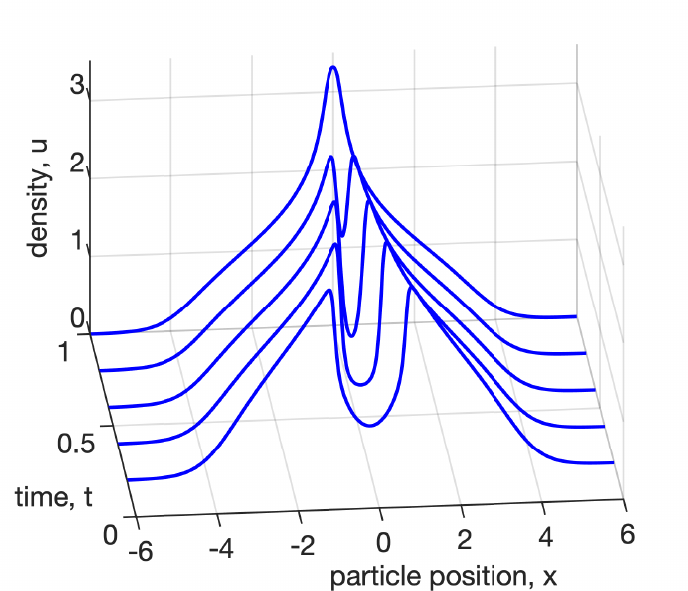}}\hspace{-.2in}
\subfigure[Time $t=1$ to $5$] %=\frac{1}{5}(2+6x+\frac{\pi}{2}\sin(2\pi x))
{\includegraphics[scale=1.1]{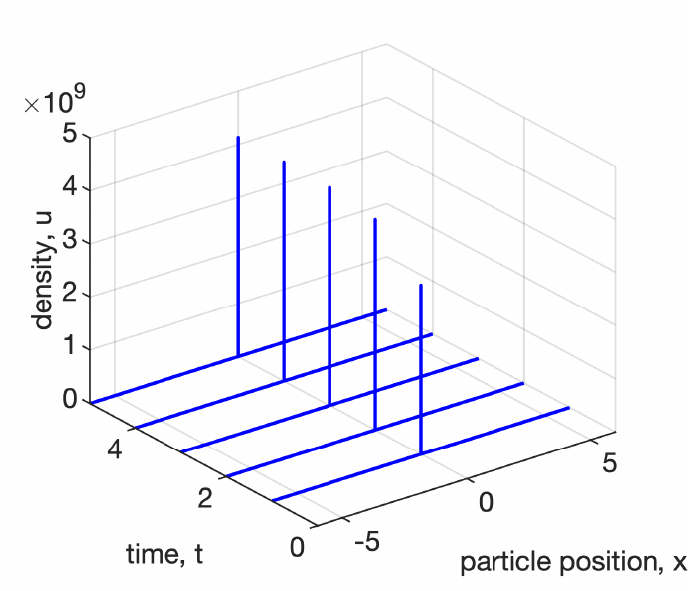}}
\caption{The evolution of density $u$ with supercritical mass $M=10$ in \textbf{Example 4}  ($h=1/1000$, $\tau=1/10000$, $m=2$).}
\label{fig:BGDensityB}
\end{figure}

%\begin{figure}
%%\begin{minipage}[t]{0.5\linewidth}
%\captionsetup{font={scriptsize}}
%\centering
%\subfigure[\scriptsize Time t=0 $\sim$ 0.15]
%{\includegraphics[width=6cm,height=5cm]{fig/BGDensity1.eps}}\hspace{0.3in}
%\subfigure[\scriptsize  Time t=0.2 $\sim$  2] %=\frac{1}{5}(2+6x+\frac{\pi}{2}\sin(2\pi x))
%{\includegraphics[width=6cm,height=5cm]{fig/BGDensity2.eps}}
%\caption{The evolution of density  for subcritical mass $M=10$ for \textbf{ Example 4}  ($h=1/10000$, $\tau=1/10000$, $m=2$)}
%\label{fig:BGDensityB}
%\end{figure}

\begin{figure}
\centering
\subfigure[Entropy]
{\includegraphics[scale=1.3]{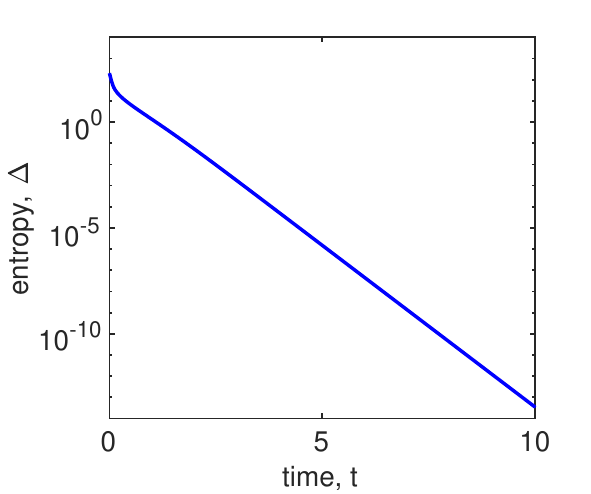}}\hspace{-.2in}
\subfigure[Total mass]
{\includegraphics[scale=1.33]{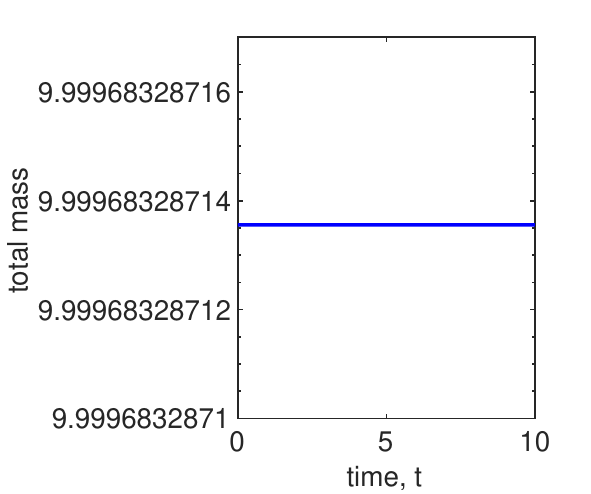}}
\caption{The evolution of entropy and total mass  with subcritical mass $M=10$ and $\beta=1$ in \textbf{Example 4}  ($h=1/1000$, $\tau=1/10000$, $m=2$).}
\label{fig:BGEnergyM=10}
\end{figure}

%\begin{figure}
%%\begin{minipage}[t]{0.5\linewidth}
%\captionsetup{font={scriptsize}}
%\centering
%\subfigure[\scriptsize Entropy]
%{\includegraphics[width=5cm,height=4cm]{fig/BGEntropy_t=20.eps}}\hspace{.3in}
%\subfigure[\scriptsize Particle] %=\frac{1}{5}(2+6x+\frac{\pi}{2}\sin(2\pi x))
%{\includegraphics[width=5cm,height=4cm]{fig/BosonGaseMassCenter_t=20.eps}}
%\caption{ The evolution of entropy and mass of the center point for subcritical mass $M=10$ and $\beta=1$  for \textbf{ Example 4}  ($h=1/10000$, $\tau=1/10000$, $m=2$)}
%\label{fig:BGEnergyB}
%\end{figure}

\begin{figure}
\centering
\subfigure[Mass at the central point]
{\includegraphics[scale=1.25]{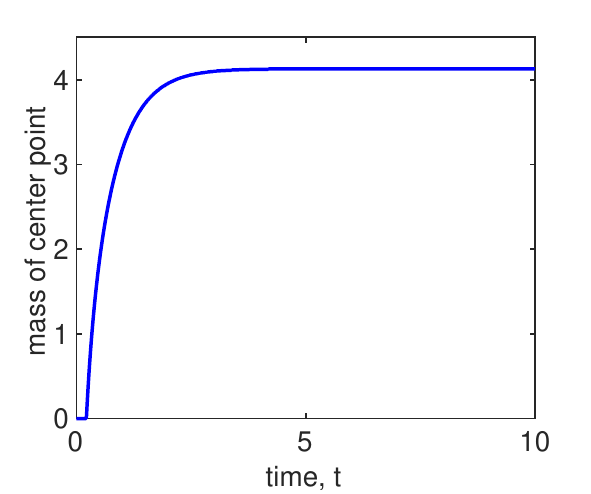}}\hspace{.2in}
\subfigure[The relationship between $\beta$ and the critical mass ]
{\includegraphics[scale=1.25]{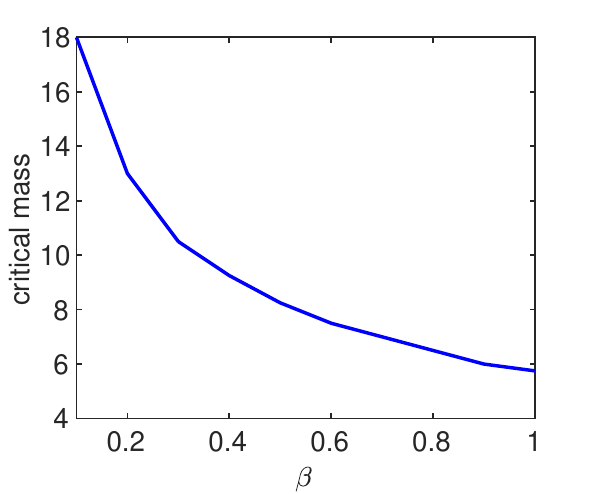}}
\caption{The evolution of mass at central point ($\beta=1$) and  the critical mass with different $\beta$  with subcritical mass $M=10$ in \textbf{Example 4}  ($h=1/1000$, $\tau=1/10000$, $m=2$).}
\label{fig:BGMass}
\end{figure}

%\begin{figure}
%\centering
%\includegraphics[scale=1.]{fig/BGMass_beta.eps}
%\caption{The relationship between $\beta$ and the critical mass in \textbf{Example 4} ($h=1/10000$, $\tau=1/10000$, $m=2$).}
%\label{fig:DWbeta}
%\end{figure}

%\begin{figure}
%%\begin{minipage}[t]{0.5\linewidth}
%\captionsetup{font={scriptsize}}
%\centering
%%\subfigure[\scriptsize The evolution of density for {\bf Example 1}]
%\includegraphics[width=5cm,height=4cm]{fig/BGMass_beta.eps}%\hspace{.3in}
%%\subfigure[\scriptsize The density with different $m$] %=\frac{1}{5}(2+6x+\frac{\pi}{2}\sin(2\pi x))
%%{\includegraphics[width=7cm,height=6cm]{fig/DWDensityM.eps}}
%\caption{The relation between $\beta$ and the critical mass  for \textbf{ Example 4} ($h=1/10000$, $\tau=1/10000$, $m=2$)}
%\label{fig:DWbeta}
%\end{figure}

%We have used the similar algorithm to solve the singularity of random genetic drift problem in \cite{C.H. Duan(2017)}.

We first consider $\beta=1$. Fig. \ref{fig:BGDensity} and \ref{fig:BGDensityB} present the density $u$ for the cases with subcritical mass $M=1$ and supercritical mass $M=10$, respectively. One can see that the solution remains bounded with subcritical mass, and the solution blows up at a finite time with supercritical mass, being consistent with the theoretical conclusion made in \cite{N.B.Abdallah(2011)} and the numerical observation in \cite{M.Bessemoulin-Chatard(2012),H.Liu(2016),Z.Sun(2018)}. One remarkable advantage of our numerical method is that the numerical solution obtained by the scheme (\ref{equ:numnum}) approximates the singularity of the scale $O(1/\varepsilon_0)$, with the small positive $\varepsilon_0$ close to the machine precision.  %Here we take $\varepsilon_0=10^{-9}$.
Fig. \ref{fig:BGEnergy} (a)  displays that the entropy decays monotonically and remains positive as time evolves for the subcritical case. Fig. \ref{fig:BGEnergy} (b)  shows the concentration process of particles towards the origin for $M=1$.  Fig. \ref{fig:BGEnergyM=10}  (a) and (b) show the evolution of  entropy  and total mass for the case with supercritical mass $M=10$.  We can see that the entropy decays exponentially to zero and the total mass remains constant. Fig. \ref{fig:BGMass}  (a)  shows  the evolution of mass at the central point for the case with supercritical mass $M = 10$ and  $\beta=1$. The increase of mass  reveals that the particles accumulate at the central  point. A saturation mass, related to the value of $\varepsilon_0$, is achieved when the concentration process balances the diffusion process.  Also, we study the relationship between $\beta$ and the critical mass in Fig. \ref{fig:BGMass} (b). As $\beta$ grows, the force exerted by the external potential $V(x)$ gets stronger and the critical mass decreases correspondingly.

\subsection{Aggregation-diffusion Models}
We now consider the Fokker--Planck equations with nonlocal interaction kernels.
%\subsubsection{Nonlinear Diffusion with a Smooth Interaction Kernel}
The Fokker--Planck equations with a smooth kernel have been studied in \cite{M.Burger(2014)}, which has proved that, under some conditions, there exists a unique steady state with a compact support. Moreover, the steady solution is a minimizer of a total energy \cite{J.Bedrossian(2011)}. In \cite{Z.Sun(2018)} and \cite{R.Bailo(2018)}, the authors verified the above results with a discontinuous Galerkin method and finite volume schemes, respectively.

{\bf Example 5: Gaussian Interaction Kernel}

In this example,  we solve the Fokker--Planck equations
\begin{equation}\label{eqiniGK}
\partial_t u=\partial_x\left[u\partial_x(\nu u^{m-1}+W*u)\right], \ x\in[-6, 6],
\end{equation}with a Gaussian kernel $W(x)=-\frac{1}{\sqrt{2\pi\sigma^2}}e^{-\frac{|x|^2}{2\sigma^2}}$,  $\sigma>0$,  and the same initial condition as in  \cite{Z.Sun(2018)}:
$$u(x,0)=\frac{1}{2\sqrt{2\pi}}\left[e^{-\frac{(x-\frac{5}{2})^2}{2}}
+e^{-\frac{(x+\frac{5}{2})^2}{2}}\right].$$
We split the Gaussian kernel $W$ as follows: $$W=W_c - W_e,$$ where $$W_c=ax^2,$$ and $$W_e=\frac{1}{\sqrt{2\pi\sigma^2}}e^{-\frac{|x|^2}{2\sigma^2}}{+}ax^2,$$ with $$a=\frac{1}{\sigma^2\sqrt{2\pi \sigma^2}}\max\left\{1,e^{-\frac{l^2}{2 \sigma^2}}\Big(1-\frac{l^2}{\sigma^2}\Big)\right\}.$$ Notice that both $W_c$ and $W_e$ are convex functions.

\begin{figure}
\centering
\subfigure[Steady States]
{\includegraphics[scale=1.3]{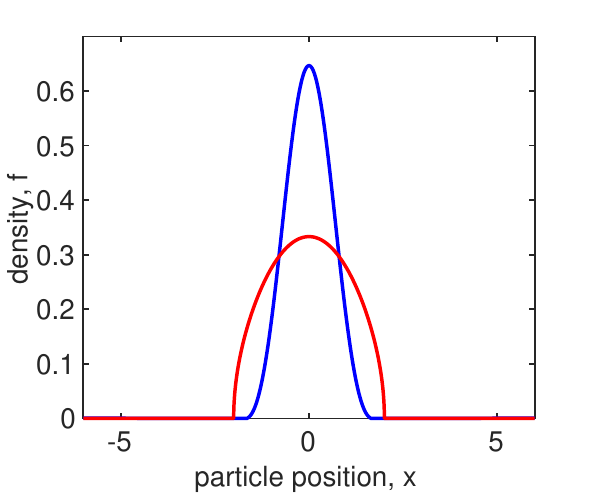}}
\subfigure[Total Energy]
{\includegraphics[scale=1.3]{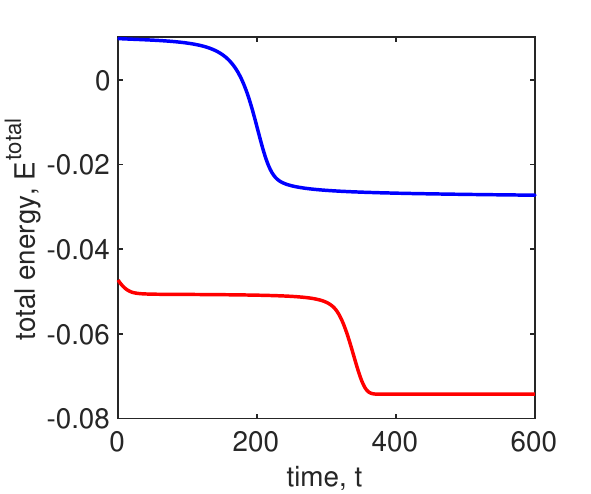}}
\caption{ Steady-state solutions and the evolution of the total energy with $m=1.5,\nu=1.28$ (the blue line) and $m=3,\nu=0.28$  (the red line) in \textbf{Example 5}  ($h=1/100$, $\tau=1/100$).}
\label{fig:GasEnergy}
\end{figure}
%   \begin{figure}
%%\begin{minipage}[t]{0.5\linewidth}
%\captionsetup{font={scriptsize}}
%\centering
%\subfigure[\scriptsize Steady State at $T=1800$]
%{\includegraphics[width=5cm,height=4cm]{fig/GaussianDensity.eps}}\hspace{.3in}
%\subfigure[\scriptsize Total Energy] %=\frac{1}{5}(2+6x+\frac{\pi}{2}\sin(2\pi x))
%{\includegraphics[width=5cm,height=4cm]{fig/GaussianEnergy.eps}}
%\caption{ Steady solution  and the evolution of total energy with $m=1.5,\nu=1.28$ and $m=3,\nu=0.28$  for \textbf{Example 5}  ($h=1/100$, $\tau=1/100$)}
%\label{fig:GasEnergy}
%\end{figure}
 \begin{figure}
\centering
\subfigure[$m=1.5$ and $\nu=0.28$]
{\includegraphics[scale=1.3]{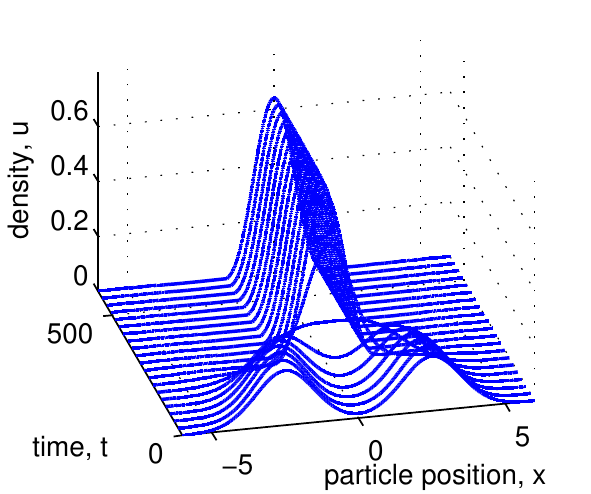}}
\subfigure[$m=3$ and $\nu=1.48$] %=\frac{1}{5}(2+6x+\frac{\pi}{2}\sin(2\pi x))
{\includegraphics[scale=1.3]{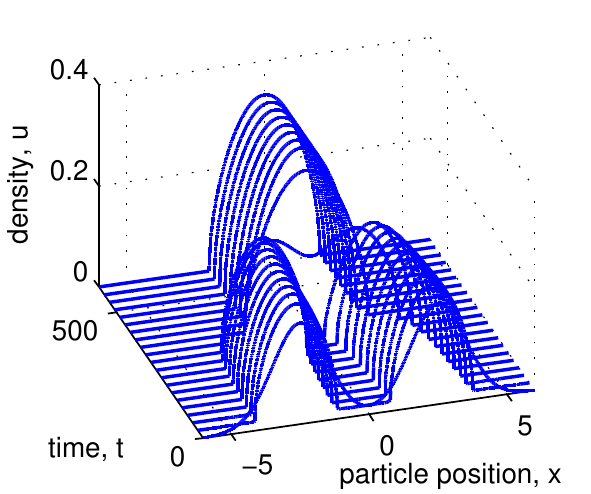}}
\caption{ Evolution of the density $u$ in \textbf{Example 5} ($h=1/100$, $\tau=1/100$)}
\label{fig:GasDenEvolution}
\end{figure}

\begin{figure}
\centering
\subfigure[$m=1.5$ and $\nu=0.28$]
{\includegraphics[scale=1.3]{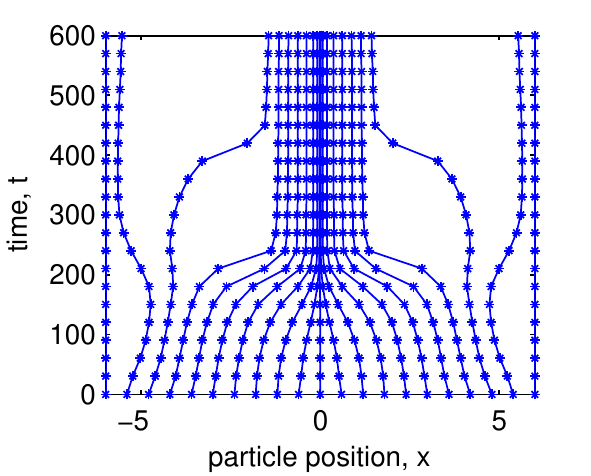}}
\subfigure[$m=3$ and $\nu=1.48$] %=\frac{1}{5}(2+6x+\frac{\pi}{2}\sin(2\pi x))
{\includegraphics[scale=1.27]{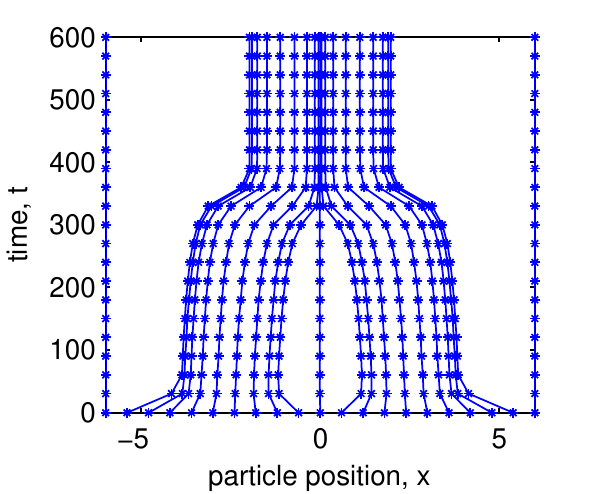}}
\caption{ Particle trajectories in \textbf{Example 5}  ($h=1/100$, $\tau=1/100$)}
\label{fig:GasParEvolution}
\end{figure}

\begin{table}[H]
\centering
\begin{tabular}{cccccc}
\hline

\hline
$h$    &$\tau$ &$\|e_h^u\|_2 $    & Order   &$\|e_h^u\|_{\infty} $  & Order  \\\hline
 1/10&1/10 &7.880e-05&&5.158e-05& \\\hline
 1/20 &1/40 &1.985e-05&1.988&1.332e-05&1.953\\ \hline
 1/40 & 1/160 &4.098e-06 &2.276&2.992e-06&2.155 \\ \hline
$h$    &$\tau$& $\|e_h^x\|_2 $     & Order   &$\|e_h^x\|_{\infty}  $     & Order\\\hline
 1/10&1/10 &6.291e-04&&6.778e-04&\\ \hline
 1/20 &1/40 &1.632e-04&1.953&1.767e-04&1.939 \\ \hline
1/40&1/160 &3.760e-05&2.155&4.177e-05&2.081 \\ \hline

  \hline
\end{tabular}
\caption{Numerical error and convergence order of the numerical trajectory $x$ and density $u$ with $m=1.5$ and $\nu=0.28$ at time $t=1$ in {\bf Example 5}.}\label{table:NonSmooth}
\end{table}

%We solve the corresponding trajectory equation \eqref{eqtra}    with the initial and boundary condition \eqref{eqtraini}-\eqref{eqtraboun}  by the fully discrete scheme \eqref{equ:numnum} with \eqref{equ:numbou} and approximate the density function $u$  in \eqref{equ:conservationL} by \eqref{numdist}-\eqref{numdist2}.

We solve the problem up to time $T=1800$, at which the system almost reaches the steady state  with $\sigma=1$.  Fig. \ref{fig:GasEnergy} (a) shows the steady-state densities with compact supports for the cases with $m=1.5, \nu=0.28$ and $m=3,\nu=1.48$.  We observe that the support gets larger as $m$ increases.  Fig. \ref{fig:GasEnergy} (b) shows that the total energy decays slowly in the first stage and then decays sharply after certain critical time, especially for the case with $m=3$. Similar results have been reported  in the work \cite{Z.Sun(2018)}, mainly due to the appearance of metastability.  Fig. \ref{fig:GasDenEvolution} (a) and (b) present the evolution of densities for $m=1.5$ and $m=3$, respectively. In comparison with the case of $m=1.5$,  the initial two smooth peaks quickly turn to two sharp bumps staying away from each other for $m=3$.  As the bumps get closer, the solution develops one support and grows up to a steady state. However, with $m=1.5$, the initial two peaks merge into one smooth peak and subsequently go to the corresponding steady state.  Fig. \ref{fig:GasParEvolution} (a) and (b) show the motion of particles for $m=1.5, \nu=0.28$ and $m=3,\nu=1.48$, respectively.  For the case of $m=1.5, \nu=0.28$, most of the particles gradually move to the origin, developing a density bump at the origin. For $m=3,\nu=1.48$, in contrast, the particles first concentrate and develop two density bumps, and subsequently move towards the origin. We also study the numerical accuracy of our schemes with the presence of nonlocal interactions. Table \ref{table:NonSmooth} shows that the  error of numerical solution $u$ and trajectory $x$ with  $m=1.5$ and $\nu=0.28$ at time $t=1$ in $\mathcal{L}^2$ and $\mathcal{L}^{\infty}$ norms. The reference ¡°exact¡± solution is obtained numerically on a rather refined mesh with $h=1/1000$ and $\tau=1/1000$. The results reveal that our numerical method is second-order accurate in space and first-order accurate in time.

\begin{rem}
The numerical simulations reveal a multi-phase convergence to equilibrium rather than a fixed-rate convergence, reminiscent of metastability.  Intermediate aggregations that depend on the initial data can quickly form, even though the final steady state is simply connected and compactly supported. These aggregations eventually merge with an arbitrarily slow convergence rate,  if the parameter $\sigma$ is small  \cite{R.Bailo(2018),Z.Sun(2018)}. %\zhou{Add a citation here.}
  \end{rem}

\section{Conclusions}\label{sec:6}
In this work, we have proposed novel structure-preserving numerical schemes, based on the Energetic Variational Approach, to robustly solve the nonlinear Fokker--Planck equations with nonlocal interactions. The trajectory equation has been obtained by using the balance between the maximal dissipation principle and least action principle.  With a convex-splitting technique, we have established numerical schemes that are uniquely solvable, with their numerical solutions satisfying the discrete energy dissipation law.  Moreover, it has been proved that the developed numerical schemes can preserve mass conservation and positivity of solutions at fully discrete level.    Numerical accuracy of second order in space and first order in time can be theoretically justified with detailed numerical analysis.

Numerical simulations have demonstrated several valuable features of the proposed schemes. In addition to the preservation of physical structures, such as positivity, conservation, discrete energy dissipation, and steady states, numerical tests have revealed that the developed numerical schemes are able to effectively and robustly solve degenerate cases of the Fokker--Planck equations with nonlocal interactions. For instance, our numerical schemes have been shown to have convergence order in degenerate cases in the presence of solutions with compact support, accurately calculate the waiting time of free boundaries without any oscillation, and approximate blow-up singularity with machine precision.

% According to the numerical simulation results, the solution could be numerically calculated  effectively without no oscillation  around the free boundary.  We   give a criterion that can compute the waiting time naturally. Another advantage is the ability to catch the blow-up singular at a machine precise. For the nonlocal kernel, the numerical results verify the analysis results.  Moreover, the  convergence order can also be obtained   for both smooth text and compacted support test numerically.

%We now discuss several issues and possible further refinements of our work.  One obvious limitation of this work is associated with the one-dimensional nature of the problem. In higher dimensions, the trajectory equation becomes a  complex nonlinear degenerate parabolic system with the determinant of the deformation gradient, i.e., $\det\frac{\partial x}{\partial X}$

 %We can simulate the finite speed and the smooth interfaces in a sole system.

We now discuss several issues and possible further refinements of our work. Our numerical method has been proved to be second-order accurate in space and first-order accurate in time. It is desirable but challenging to develop second-order accurate temporal discretization that is able to preserve unconditional energy dissipation in the discrete sense. One promising idea is to employ the BDF discretization with an artificial Douglas-Dupont regularization term, which is added to ensure the energy dissipation \cite{Chen_JCP2019}. Another improvement is to develop numerical methods and analysis to address a non-smooth (even singular) interaction kernel $W(\cdot)$. One possible strategy is to split the kernel into a smooth part and a non-smooth part that can be treated analytically.  Finally, one limitation of this work is associated with the one-dimensional nature of the problem. In higher dimensions, the trajectory equation will be a very complicated nonlinear parabolic system with the Jacobian of the flow map in the denominator.  The development of numerical schemes with structure-preserving properties for higher dimensional cases deserves further investigation.

%Another limitation of this work is associated with the one-dimensional nature of the problem. In two or higher dimension, the trajectory equation will be a  complex nonlinear degenerate parabolic equation system with the determinant of the deformation gradient, i.e., $\det\frac{\partial x}{\partial X}$.   We are  finding a suitable numerical method which can satisfy the discrete energy dissipation law and give the convergence analysis. Solving for multi-dimensional these Fokker--Planck type equations by this energetic method will be left to our future works.

%Another work is establishing an algorithm for the non-smooth nonlocal term $W$, especially for the singular kernel. We have established a similar numerical scheme for the non-smooth $W$ and dealt with the singular part by a splitting technical. Moreover, the unique solvability can be also proven. Relevant results will be presented elsewhere.

\section*{Acknowledgments}
The authors would like to thank anonymous reviewers for their helpful suggestions which lead to improvement of the work. C. Duan was supported in part by NSFC under the grant 11901109. C. Liu was partially supported by the United States-Israel Binational Science Foundation (BSF) \# 2024246, and NSF grants DMS-1216938 and DMS-1418689. W. Chen was supported by the National Science Foundation of China (11671098) and partially supported by Shanghai science and technology research program (19JC1420101). W. Chen also thanks Institute of Scientific computation and Financial Data Analysis, Shanghai University of Finance and Economics for the support during his visit. X. Yue was partially supported by NSFC under the grant 11971342. S. Zhou was supported by the grants NSFC 21773165, Young Elite Scientist Sponsorship Program by Jiangsu Association For Science and Technology, and National Key R\&D Program of China  (No. 2018YFB0204404).

%%-----------------------------
%%      your bibliography
%%-----------------------------
%\footnotesize
\bibliographystyle{unsrt}
%\bibliography{EnVarFP}

\begin{thebibliography}{99}

\bibitem{N.B.Abdallah(2011)}
N. B. Abdallah, I. M. Gamba, G. Toscani, On the minimization problem of sub-linear convex functionals, Kinet. Relat. Models 4(4) (2011) 857-871.

\bibitem{D.G.Aronson(1983)}
D. G. Aronson, L. A. Caffarelli, S. Kamin, How an initially stationary interface begins to move in porous medium flow, SIAM J. Math. Anal. 14(4) (1983) 639-658.


\bibitem{R.Bailo(2018)}
R. Bailo, J. A. Carrillo,  J. Hu, Fully discrete positivity-preserving and energy-decaying schemes for aggregation-diffusion equations with a gradient flow structure, arXiv preprint (2018) arXiv:1811.11502.

\bibitem{J.Barre(2017)}
J. Barr\'{e}, P. Degond, E. Zatorska, Kinetic theory of particle interactions mediated by dynamical networks, Multiscale Model. Simul. 15(3) (2017) 1294-1323.

\bibitem{J.Bedrossian(2011)}
J. Bedrossian, Global minimizers for free energies of subcritical aggregation equations with
degenerate diffusion, Appl. Math. Lett. 24(11) (2011) 1927-1932.

\bibitem{D.Benedetto(1998)}
D. Benedetto, E. Caglioti, J. A. Carrillo, M. Pulvirenti, A non-Maxwellian steady distribution for one-dimensional granular media, J. Stat. Phys. 91 (1998) 979-990.

%\bibitem{D.Benedetto(1997)}
%D. Benedetto, E. Caglioti, M. Pulvirenti, A kinetic equation for granular media, RAIRO. Mod Å¡Å l. Math. Anal. Num\'{e}r. 31(5) (1997) 615-641.

\bibitem{M.Bessemoulin-Chatard(2012)}
M. Bessemoulin-Chatard, F. Filbet, A finite volume scheme for nonlinear degenerate parabolic equations, SIAM J. Sci. Comput. 34(5) (2012) B559-B583.

\bibitem{S.Boscarino(2016)}
S. Boscarino, F. Filbet, G. Russo, High Order Semi-implicit Schemes for Time Dependent Partial Differential Equations, J. Sci. Comput. 68 (2016) 975-1001.

\bibitem{C.Buet(2004)}
C. Buet, S. Cordier, and V. Dos Santos, A conservative and entropy scheme for a simplified model of granular media, Transport  Theor. Stat. 33(2) (2004) 125-155.



\bibitem{M.Burger(2014)}
M. Burger, R. Fetecau,  Y. Huang, Stationary states and asymptotic behavior of aggregation
models with nonlinear local repulsion, SIAM J. Appl. Dyn. Syst. 13(1) (2014) 397-424.

\bibitem{M.Burger(2007)}
M. Burger, V. Capasso, D. Morale, On an aggregation model with long and short range interactions, Nonlinear Anal. Real World Appl. 8(3) (2007) 939-958.

%\bibitem{M.Burger(2010)}
%M. Burger, J. A. Carrillo, M. T. Wolfram, A mixed finite element method for nonlinear diffusion equations, Kinet. Relat. Models 3(1) (2010) 59-83.

\bibitem{J.Carrillo(2001)}
J. A. Carrillo, A. J$\ddot{u}$ngel, P. A. Markowich, G. Toscani, A. Unterreiter, Entropy dissipation methods for degenerate parabolic problems and generalized Sobolev inequalities, Monatshefte Math. 133(1) (2001) 1-82.

\bibitem{J.Carrillo(2000)}
J. A. Carrillo, G. Toscani, Asymptotic $L^1$-decay of solutions of the porous medium equation to self-similarity, Indiana Univ. Math. J. 49(1) (2000) 113-142.

\bibitem{J.Carrillo(2003)}
J. A. Carrillo, R. J. McCann, C. Villani, Kinetic equilibration rates for granular media and related equations: entropy dissipation and mass transportation estimates, Rev. Mat. Iberoam. 19(3) (2003) 971-1018.

\bibitem{J.Carrillo(2014)}
J. A. Carrillo, Y. Huang, S. Martin, Explicit
flock solutions for Quasi-Morse potentials, Eur. J. Appl. Math., 25(5) (2014) 553-578.


\bibitem{J.Carrillo(2010)}
J. A. Carrillo, M. Fornasier, G. Toscani,  F. Vecil, Particle, kinetic, and hydrodynamic models
of swarming. In Mathematical modeling of collective behavior in socio-economic and life sciences,
Model. Simul. Sci. Eng. Technol., pages 297-336. Birkh\"{a}user Boston, Inc., Boston, MA, 2010.


\bibitem{J.Carrillo(2015)}
J. A. Carrillo, A. Chertock,  Y. Huang, A finite-volume method for nonlinear nonlocal equations with a gradient flow structure, Commun. Comput. Phys.  17 (2015) 233-258.

\bibitem{J.Carrillo(2016)}
J. A. Carrillo, H. Ranetbauer, M. T. Wolfram, Numerical simulation of nonlinear continuity equations by evolving diffeomorphisms, J. Comput. Phys. 327 (2016) 186-202.

\bibitem{J.Carrillo(2017)}
J. A. Carrillo, Y. Huang, F. S. Patacchini, G. Wolansky, Numerical study of a particle method for gradient flows, Kinet. Relat. Models 10(3) (2017) 613-641.

\bibitem{J.Carrillo(2018)}
J. A. Carrillo, K. Craig, Y. Yao, Aggregation-diffusion equations: dynamics, asymptotics, and singular limits, arXiv preprint arXiv:1810.03634, 2018.

\bibitem{Chen_JCP2019}
W. Chen, C. Wang, X. Wang, and S. M. Wise, Positivity-preserving, energy stable numerical schemes for the Cahn--Hilliard equation with logarithmic potential, J. Comput. Phys.:X 3  (2019) 100031.
%\bibitem{J.A.Carrillo(2018)}
%J. A. Carrillo, A. Colombi, M. Scianna, Adhesion and volume constraints via nonlocal interactions
%determine cell organisation and migration profiles, J. Theoret. Biol.  445 (2018) 75-91.

%\bibitem{K.Craig(2016)}
%K. Craig, A. Bertozzi, A blob method for the aggregation equation, Math. Comp. 85(300) (2016) 1681-1717.

\bibitem{Q.Du(2009)}
Q. Du,  C. Liu, R. Ryham,  X. Wang, Energetic variational approaches in modeling vesicle and fluid interactions, Phys. D 238 (2009) 923-930.

\bibitem{C.H.Duan(2017)}
C. Duan, C. Liu, C. Wang,  X.  Yue, Numerical complete solution for random genetic drift by Energetic Variational approach, ESAIM: Math.  Model.  Num. 53(2) (2019)  615-634.

\bibitem{C.H.Duan(2018)}
C. Duan, C. Liu, C. Wang, X.  Yue, Numerical methods for Porous Medium Equation by an Energetic Variational Approach, J. Comput. Phys.   385  (2019) 13-32.

\bibitem{C.H.Duan(2019)}
C. Duan, C. Liu, C. Wang, X.  Yue,  Convergence Analysis of a Numerical Scheme for the Porous Medium Equation by an Energetic Variational Approach, Numer. Math. Theor. Meth. Appl. 13 (2020).

\bibitem{C.H.Duan(2020)}
C. Duan, W. Chen, C. Liu, C. Wang, X. Yue, A second order accurate numerical scheme for the porous medium equation by an energetic variational approach, arXiv preprint  arXiv:2006.12354 (2020).


\bibitem{W.E(1995)}
W. E, J. G. Liu, Projection method I: convergence and numerical boundary layers, SIAM J. Numer. Anal. 32 (1995) 1017-1057.


\bibitem{B.Eisenberg(2010)}
B. Eisenberg, Y. K. Hyon,  C. Liu, Energy variational analysis of ions in water and channels: Field theory for primitive models of complex ionic fluids,  J. Chem. Phys. 133(10)  (2010) 104.

\bibitem{D.J.Eyre(1998)}
D. J. Eyre, Unconditionally gradient stable time marching the Cahn-Hilliard equation, in \emph{MRS Proceedings}, Cambridge Univ. Press 529 (1998)  39.


\bibitem{Y.Hyon(2010)}
Y. Hyon, D. Y.  Kwak, C. Liu, Energetic variational approach in complex fluids: maximum dissipation principle, Discrete Contin. Dyn. Syst.  26(4) (2010) 1291-1304.

\bibitem{Koba(2017)}
 H. Koba, C. Liu, Y. Giga, Energetic variational approaches for incompressible fluid systems on an evolving surface, Quart. Appl. Math. 75 (2017)  359-389.

\bibitem{T.Kolokolnikov(2013)}
T. Kolokolnikov, J. A. Carrillo, A. Bertozzi, R. Fetecau, and M. Lewis. Emergent behaviour in
multi-particle systems with non-local interactions [Editorial], Phys. D 260 (2013) 1-4.

\bibitem{C.Liu(2003)}
C. Liu, J. Shen, A phase field model for the mixture of two incompressible fluids and its approximation by a Fourier-spectral method, Phys. D  179(3-4) (2003)  211-228.

%\bibitem{H. Liu(2016)}
%H. Liu, Z. Wang, An entropy satisfying discontinuous Galerkin method for nonlinear Fokker--Planck equations, J. Sci. Comput.  68(3) (2016) 1217-1240.

 \bibitem{C.Liu(2017)}
C. Liu,  H. Wu,  An energetic variational approach for the Cahn-Hilliard equation with dynamic boundary conditions, Arch. Ration. Mech. Anal.  233(1)  (2019) 167-247.

 %\bibitem{C.Liu(2019)}
% C. Liu,  Y. Wang, On Lagrangian schemes for the multidimensional porous medium equations by a discrete energetic variational approach,  arXiv preprint (2019) arXiv:1905.12225 .

 \bibitem{H.Liu(2016)}
H. Liu, Z. Wang, An entropy satisfying discontinuous Galerkin method for nonlinear Fokker--Planck equations, J. Sci. Comput. 68(3) (2016) 1217-1240.

\bibitem{H.Liu(2017)}
H. Liu, Z. Wang, A free energy satisfying discontinuous Galerkin method for one-dimensional Poisson-Nernst-Planck systems, J. Comput. Phys. 328 (2017) 413-437.


 \bibitem{P.M.Lushnikov(2008)}
P. M. Lushnikov, N. Chen,   M. Alber, Macroscopic dynamics of biological cells interacting via
chemotaxis and direct contact,  Phys. Rev. E 78 (2008) 061904.

 \bibitem{Y.Nesterov(1994)}
Y. Nesterov, A. Nemirovskii,  Interior-point polynomial algorithms in convex programming, SIAM, 1994.

\bibitem{L.Onsager(1931)}
L. Onsager, Reciprocal relations in irreversible processes, Phys. Rev. II. Ser. 38 (1931) 2265-2279.

\bibitem{L.Onsager1(1931)}
L. Onsager, Reciprocal relations in irreversible processes, Phys. Rev. I. 37(4) (1931) 405.

\bibitem{L.Pareschi(2007)}
L. Pareschi,  M. Zanella, Structure Preserving Schemes for Nonlinear Fokker-Planck Equations and Applications, J. Sci. Comput. 73(3) (2017) 1575-1600.

\bibitem{Y.Qian(2019)}
Y. Qian, Z. Wang, S. Zhou, A conservative, free energy dissipating, and positivity preserving finite difference scheme for multi-dimensional nonlocal Fokker--Planck equation, J. Comput. Phys.  386 (1) (2019)  22-36.

\bibitem{LiuWangJCP(2017)}
H. Liu, Z. Wang, A free energy satisfying discontinuous Galerkin method for one-dimensional {Poisson--Nernst--Planck} systems, J. Comput. Phys. 328 (2017) 413-437.

\bibitem{LiuWangJCP(2014)}
H. Liu, Z. Wang, A free energy satisfying finite difference method for {Poisson--Nernst--Planck}  equations, J. Comput. Phys. 268 (2014) 363-376.

\bibitem{MettiXuLiuJCP(2016)}
M. S. Metti, J. Xu, C. Liu, Energetically stable discretizations for charge transport and electrokinetic models, J. Comput. Phys. 306 (2016) 1-18.


\bibitem{DingWangZhou_JCP19}
J. Ding, Z. Wang, S. Zhou, Positivity preserving finite difference methods for {Poisson--Nernst--Planck} equations with steric interactions: Application to slit-shaped nanopore conductance, J. Comput. Phys.  397  (2019)  108864.

\bibitem{J.W.Strutt(1873)}
J. W. Strutt, Some general theorems relating to vibrations, P. Lond. Math. Soc. IV (1873) 357-368.

\bibitem{Z.Sun(2018)}
Z. Sun, J. A. Carrillo, C.-W. Shu, A discontinuous Galerkin method for nonlinear parabolic equations and gradient flow problems with interaction potentials,  J. Comput. Phys. 352 (2018) 76-104.

\bibitem{C.M.Topaz(2006)}
C. M. Topaz, A. L. Bertozzi,   M. A. Lewis, A nonlocal continuum model for biological aggregation,
Bull. Math. Biol.  68 (2006) 1601-1623.



\bibitem{J.L.Vazquez(2007)}
J. L. V\'azquez, The Porous Medium Equation, Oxford University Press, Oxford, 2007.

\bibitem{C.Villani(2003)}
C. Villani, Topics in Optimal Transportation, American Mathematical Society, 2003.

\bibitem{C.Wang(2000)}
C. Wang, J. G. Liu, Convergence of gauge method for incompressible flow, Math. Comp.  69 (2000) 1385-1407.

\end{thebibliography}

\end{document}